\documentclass[11pt,a4paper]{article}
\usepackage[a4paper,margin=29mm]{geometry}
\usepackage[T1]{fontenc}
\usepackage[utf8]{inputenc}
\usepackage{mathpazo,microtype}
\usepackage{amsmath,amssymb,amsthm,mathtools,enumitem}
\usepackage{xcolor}
\usepackage{hyperref}
\hypersetup{colorlinks=true,linkcolor=blue!45!black,citecolor=blue!45!black,
 urlcolor=blue!45!black,pdftitle={Morse theory and a degree formula for the resonant singular Paneitz equation on four-manifolds},
 pdfauthor={Mohameden Ahmedou, Lina Wu, Lei Zhang},
 pdfsubject={Critical points at infinity, Morse inequalities, and resonant degree}}
\numberwithin{equation}{section}
\newtheorem{theorem}{Theorem}[section]
\newtheorem{lemma}[theorem]{Lemma}
\newtheorem{proposition}[theorem]{Proposition}
\newtheorem{corollary}[theorem]{Corollary}
\newcommand{\R}{\mathbb R}
\newcommand{\FF}{\mathbb F_2}
\newcommand{\F}{\mathcal F}
\newcommand{\B}{\mathcal B}
\newcommand{\cH}{\mathcal H}
\newcommand{\dV}{\,dV_g}
\newcommand{\Crit}{\operatorname{Crit}}
\newcommand{\Ent}{\operatorname{Ent}}
\newcommand{\ind}{\operatorname{ind}}
\newcommand{\Conf}{\operatorname{Conf}}
\newcommand{\Cov}{\operatorname{Cov}}
\newcommand{\fp}{\mathrm{fp}}
\newcommand{\one}{\mathbf 1}
\title{\bfseries Morse theory and a degree formula for the resonant singular Paneitz equation on four-manifolds}
\author{Mohameden Ahmedou \and Lina Wu \and Lei Zhang}
\newcommand{\shorttitle}{Morse theory for the singular Paneitz equation}
\date{20 November 2026}
\begin{document}
\maketitle
\pagestyle{myheadings}
\markboth{\shorttitle}{\shorttitle}
\begin{abstract}
We develop Morse theory for the singular Paneitz equation at resonant
mass on closed four-manifolds with conical orders between minus one
and zero and Paneitz kernel consisting of the constants. Under nondegeneracy of the reduced
interactions and nonvanishing of the leading mass-defect coefficients,
we classify the critical points at infinity by the sign of these
coefficients. We compute their indices and their contributions to
the relative homology of sublevels. The topology of very low sublevels
is determined by weighted barycenters of strictly smaller total weight
and by the negative Paneitz eigenspace. These results yield strong
Morse inequalities and a Leray--Schauder degree formula at resonance,
which also applies when genuine solutions are degenerate. On the round
four-sphere with two conical points, the Morse inequalities give
existence on an open set of data for which the degree vanishes.
\end{abstract}


\tableofcontents
\clearpage
\section{The equation, the resonant degree and the Morse problem}
At a resonant mass, bounded-energy descent for the singular
Paneitz functional can lead to concentration at ordinary points
and conical marks. We identify the concentrating configurations
that persist under a descending pseudogradient, and compute the
change of sublevel topology they cause. A reduction on slices of
fixed logarithmic height links the sign of the height derivative
to the mass-defect expansion in
\emph{The singular Paneitz equation on four-manifolds: refined
blow-up analysis and degree jumps}~\cite{AWZBlowup}.
This link allows us to treat distinct representations of the
same resonant mass together and to obtain a degree formula at
resonance from the Morse relations.

\subsection{The singular equation}
Let $(M,g)$ be a smooth closed connected four-manifold of volume $V$, and assume that $\ker P_g=\R$. Let $q_1,\ldots,q_N$ be distinct points with orders $-1<\gamma_j<0$, set $\Sigma=\{q_1,\ldots,q_N\}$ and $D=\sum_j\gamma_j\delta_{q_j}$, and fix $h\in C^\infty(M)$ positive. The desingularized equation with parameter $\kappa>0$ is
\begin{equation}\label{normalized}
 P_gu=2\kappa(\mu_u-V^{-1}dV_g),\qquad
 \mu_u=\frac{K_De^{4u}\,dV_g}{\int_M K_De^{4u}\dV},
 \qquad \int_M u\dV=0 .
\end{equation}
We seek $u$ in $\cH=\{u\in H^2(M):\int_Mu\dV=0\}$ and identify functions with their volume densities in distributional equations. The weight $K_D$ in \eqref{weight} has order $d_g(x,q_j)^{4\gamma_j}$ at $q_j$.

Put $\mu_j=1+\gamma_j$, $w(I)=\sum_{j\in I}\mu_j$ with $w(\varnothing)=0$, and $\rho=\kappa/(8\pi^2)$. Concentration for exact solutions can occur only at a positive element of $8\pi^2\mathcal W$, where
\[
 \mathcal W=\{p+\textstyle\sum_{j\in I}\mu_j:
                 p\in\mathbb N_0,\ I\subset\{1,\ldots,N\}\}.
\]
We fix such a wall $\kappa_*=8\pi^2\rho_*$. A representation of the wall is a pair $(p,I)$ with
\begin{equation}\label{wall}
 \rho_*=p+w(I),\qquad n=p+|I|.
\end{equation}
It specifies $p$ distinct ordinary centers in $X=M\setminus\Sigma$ and the prescribed marked centers in $I$. The atoms of the source $2\kappa\mu_u$ have masses $16\pi^2$ and $16\pi^2\mu_j$, respectively. These masses follow from the quantization theorem
\cite[Theorem~5.1]{JFA} and the localized dilation identity
\cite[Lemma~10.1]{AWZBlowup}, which makes each ordinary atom simple. Distinct marked subsets giving the same sum are retained individually. The ordinary centers are unordered.

\subsection{Green function and wall data}\label{sub:geometric}
Our conventions are $\Delta_g=\operatorname{div}_g\nabla$, $\delta=d^*$, and
\[
 P_g=\Delta_g^2+\delta\big((\tfrac23 R_g g-2\operatorname{Ric}_g)d\,\cdot\,\big),
 \qquad Q_g=-\tfrac1{12}(\Delta_gR_g-R_g^2+3|\operatorname{Ric}_g|^2).
\]
Thus $P_gv+2Q_g=2Q_{e^{2v}g}e^{4v}$ for smooth $v$. Put
$\kappa_P=\int_M Q_g\dV$ and $\kappa_D=\kappa_P+4\pi^2\sum_j\gamma_j$.
Let $G_0(a,\cdot)$ be the mean-zero Green function satisfying
$P_g G_0(a,\cdot)=\delta_a-V^{-1}dV_g$, and let $f$ be the mean-zero solution of
$P_gf=2\kappa_P/V-2Q_g$. Define
\begin{equation}\label{weight}
 K_D=h\exp\!\left(4f-32\pi^2\sum_j\gamma_jG_0(q_j,\cdot)\right),\qquad
 G=16\pi^2G_0 .
\end{equation}
Near its diagonal $G(a,x)=2\log(1/d_g(a,x))+H_g(a,x)$. Near $q_j$ we write
$K_D(x)=d_g(x,q_j)^{4\gamma_j}K_j(x)$ with $K_j>0$. The joint regularity of the recentered Green remainder needed for its second jets is established in \cite[\S2]{AWZBlowup}.

For a representation $(p,I)$ let $A$ denote the unordered set of ordinary centers and enumerate all its ordinary and active marked centers as $z_s$, $1\le s\le n$. Write $\mu_s$ for the corresponding local weight (one for an ordinary center), $M_s=16\pi^2\mu_s$, and $\widehat K_s=K_D$ at an ordinary center, $\widehat K_s=K_j$ at $q_j$. With the centers held fixed, put
\begin{align}
 \Phi_{A,I}(x)&=\sum_s\mu_s G(z_s,x),\nonumber\\
 \Theta_s(x)&=\widehat K_s(x)\exp\!\left(4\mu_s H_g(z_s,x)
                           +4\sum_{t\ne s}\mu_tG(z_t,x)\right),\\
 \F_{p,I}(A)&=\sum_s\mu_s\left[\log\widehat K_s(z_s)+2\mu_sH_g(z_s,z_s)
                                +2\sum_{t\ne s}\mu_tG(z_s,z_t)\right].\label{interact}
\end{align}
Near $z_s$, $K_De^{4\Phi_{A,I}}=d_g(x,z_s)^{-4-4\mu_s}\Theta_s(x)$. For an ordinary center $a_i$, symmetry of $G$ gives
$\nabla_{a_i}\F_{p,I}=\nabla_x\log\Theta_i(x)|_{x=a_i}$. When $p=0$ the value of \eqref{interact} is a constant on the one-point configuration space.

For each weight let $U_s$ be the radial, finite-volume normal solution of $\Delta^2 U_s=6|y|^{4(\mu_s-1)}e^{4U_s}$ on $\R^4$, with $U_s(0)=0$ and center zero \cite[\S3]{AWZBlowup}. For an ordinary core $U_s(y)=-\log(1+|y|^2/4)$. For each core define the far-field constant $c_s$ and radial probability $\nu_s$ by
\[
 U_s(y)=-2\mu_s\log |y|+c_s+o(1),\qquad
 d\nu_s=\frac{6}{M_s}|y|^{4(\mu_s-1)}e^{4U_s(y)}\,dy .
\]
Put $m=\max_s\mu_s$, $a_s=1/(4\mu_s)$, $h_s=\Theta_s(z_s)e^{-4c_s}/3$, and, when $\mu_s>1/2$, $m_{2,s}=\int_{\R^4}|y|^2\,d\nu_s$. Define
\[
 \mathcal L_s=\Delta_g\log\Theta_s(z_s)-\frac{2(2-\mu_s)}9R_g(z_s)
       +\begin{cases}
          |\nabla\log\Theta_s(z_s)|^2/[2(1-\mu_s)],&\mu_s<1,\\
          0,&\mu_s=1 .
        \end{cases}
\]
All derivatives here are in the observation variable. In the normalization of \cite[\S10.2]{AWZBlowup}, the coefficient and its positive rate are
\begin{align}
 C(A)&=
 \begin{cases}
 \displaystyle\frac{1}{8\kappa_*}\sum_{\mu_s=m}
       a_s M_s m_{2,s}h_s^{1/(2m)}\mathcal L_s, &m>1/2,\\[4pt]
 \displaystyle\frac{|\mathbb S^3|}{4\kappa_*}\sum_{\mu_s=1/2}
           \frac{a_s}{\mu_s}\Theta_s(z_s)\mathcal L_s, &m=1/2,\\[4pt]
 D_I^{(4)}/\kappa_*, &m<1/2,
 \end{cases}\label{allC}\\
 S_m(T)&=
 \begin{cases}e^{-2T/m},&m>1/2,\\ T e^{-4T},&m=1/2,\\e^{-4T},&m<1/2.\end{cases}\label{order}
\end{align}
In the third case $p=0$, and the finite part is
\begin{equation}\label{finitepart}
 D_I^{(4)}=\lim_{r\downarrow0}\left[
 \int_{M\setminus\cup_{j\in I} B_r(q_j)} K_De^{4\Phi_{\varnothing,I}}\dV
  -\sum_{j\in I}\frac{|\mathbb S^3|\Theta_j(q_j)}{4\mu_j}r^{-4\mu_j}\right].
\end{equation}

The value associated with a configuration is
\begin{equation}\label{limitenergy}
 c(A)=\sum_s M_s\left[\int U_s\,d\nu_s+c_s+\tfrac12\log3\right]
      -8\pi^2\F_{p,I}(A)-\kappa_*\log(\kappa_*/Z_0),
\end{equation}
where $Z_0=\int_M K_D\dV$.

\subsection{Critical configurations and the relative topology}
For each $(p,I)$, let $A$ range over critical points of the interaction $\F_{p,I}$ in \eqref{interact} on $\Conf_p(X)/\mathfrak S_p$, where $\Conf_p(X)$ consists of the ordered $p$-tuples of distinct points of $X$. Put
\[
 k_A=\ind(\F_{p,I},A),\qquad j_A=(-1)^{n+k_A}.
\]
For $p=0$ the configuration space is a point and $k_A=0$. The scalar $C(A)$ in \eqref{allC} is normalized by the relative defect $(\kappa-\kappa_*)/\kappa_*$. We impose the following condition on $h$ at the chosen wall.
\begin{equation*}\tag{W}\label{W}
 \F_{p,I}\text{ is Morse for every representation, and }
 C(A)\ne0\text{ for every critical configuration.}
\end{equation*}
Write $dm_D=K_D\,dV_g/Z_0$ and
\[
 J_\kappa(u)=Q(u)-\kappa\log Z(u),\qquad
 Q(u)=\langle P_gu,u\rangle,\qquad Z(u)=\int_M e^{4u}dm_D .
\]
Let $E^-$ and $E^+$ be the negative and positive spectral spaces
on $\cH$, set $k=\dim E^-$, and use the reduced functional
$F_\kappa(w)=\max_{z\in E^-}J_\kappa(w+z)$ on $E^+$. The maximum
is unique. Critical points correspond to those of $J_\kappa$
and the full Morse index exceeds the reduced index by $k$
(Section~\ref{sub:spectral}). For a functional $g$ write $g^a=\{g\le a\}$.

A critical point at infinity will mean a class of nonprecompact
forward orbits of the descending field below with energy bounded
from below. Two such orbits are equivalent when their limiting
weighted supports agree.

\begin{theorem}[Classification of the critical points at infinity]
\label{thm:cpi}
Under \eqref{W}, $J_{\kappa_*}$ has a bounded locally Lipschitz
descending pseudogradient with complete forward flow. A forward orbit
with energy bounded from below either has precompact image and
approaches the genuine critical set or converges in its bubble
coordinates to a critical point at infinity. At the latter end the
probability measures converge weakly to
$\sum_s M_s\delta_{z_s}/(2\kappa_*)$. The critical points
at infinity are in bijection with the unordered critical
configurations $A$ for which $C(A)<0$. Their limiting energies
are $c(A)$ in \eqref{limitenergy}.
\end{theorem}

At a genuine solution let $i(u)$ be the Morse index of $F_{\kappa_*}$ at its positive spectral component. Nondegeneracy means that the mean-zero Hessian of $J_{\kappa_*}$ has trivial kernel. For $C(A)<0$ set
\begin{equation}
 i_\infty(A)=4p-k_A+n-1.
\end{equation}
The term $4p-k_A$ counts the negative directions of $-\F_{p,I}$ at the
free centers. The other $n-1$ come from relative core masses. The full index at infinity is $k+i_\infty(A)$, and
$(-1)^{i_\infty(A)}=-j_A$.

Set $J=J_{\kappa_*}$ and use coefficients in
$\FF=\mathbb Z/2\mathbb Z$ for the local and global Morse relations.
Let
\[
 \mathcal C=\{J(u):DJ(u)=0\}\ \cup\
          \{c(A): A\in\Crit\F_{p,I},\ C(A)<0\}.
\]
This set is finite when all genuine solutions are nondegenerate.

\begin{theorem}[Change across a level]\label{change56}
Assume \eqref{W} and that all genuine solutions are nondegenerate.
If $[a,b]\cap\mathcal C=\varnothing$, $J^a$ is a strong deformation retract of $J^b$. If $c_*\in\mathcal C$, take $\varepsilon>0$ sufficiently small that the band $[c_*-\varepsilon,c_*+\varepsilon]$ meets no other element of $\mathcal C$. The relative homology of its endpoint sublevels in degree $q$ is
\begin{equation}\label{full-local-sum}
 \bigoplus_{\substack{DJ(u)=0,\ J(u)=c_*\\ k+i(u)=q}}\FF
 \quad\oplus\quad
 \bigoplus_{\substack{C(A)<0,\ c(A)=c_*\\
                   k+4p-k_A+n-1=q}}\FF .
\end{equation}
All ordinary centers are counted on the unordered configuration space. Configurations with $C>0$ give no summands.
\end{theorem}
\subsection{Endpoint topology, Morse inequalities and resonant degree}
Write $B=\B_{<\rho_*}$ for the space of probability barycenters on
$M$ whose distinct support points have total cost strictly less
than $\rho_*$. An ordinary point has cost one and a marked point
$q_j$ has cost $\mu_j$. The topology is the weak topology of
measures. The space is empty exactly at the first positive wall.

\begin{theorem}[Topology of the low sublevels]\label{main}
Let $\kappa_*>0$ be any wall \eqref{wall} and put $k=\dim E^-$, $B=\B_{<\rho_*}$.
There is $R_0$ such that for every $R\ge R_0$ the following hold.
\begin{enumerate}[label=\textup{(\roman*)},leftmargin=2em]
\item If $B\ne\varnothing$, then $F_{\kappa_*}^{-R}\simeq B$. If $B=\varnothing$, then $F_{\kappa_*}^{-R}=\varnothing$.
\item If $k=0$, $J_{\kappa_*}^{-R}=F_{\kappa_*}^{-R}$. If $k>0$, then
\[
 J_{\kappa_*}^{-R}\simeq
 \begin{cases}
 B*S^{k-1}, & B\ne\varnothing,\\
 S^{k-1}, & B=\varnothing .
 \end{cases}
\]
\end{enumerate}
The inclusions between any two sufficiently negative sublevels, reduced or full, are homotopy equivalences.
\end{theorem}
For large $R$, the relative Poincar\'e polynomial of
$(E^+,F_{\kappa_*}^{-R})$ is
\begin{equation}\label{pol}
 P_*(t)=
 \begin{cases}
 1,&B=\varnothing,\\
 t\displaystyle\sum_{r\ge0}\dim\widetilde H_r(B;\FF)\,t^r,&B\ne\varnothing .
 \end{cases}
\end{equation}
In the full space the index shift is $k$. In particular the
reduced zeroth homology of a disconnected $B$ contributes in
degree one of $P_*$.

\begin{theorem}[Strong Morse inequalities]\label{graded-morse}
Under \eqref{W}, if every genuine solution is nondegenerate, then
\begin{equation}\label{morse}
 \sum_u t^{i(u)}+\sum_{C(A)<0}t^{i_\infty(A)}
       =P_*(t)+(1+t)\mathcal Q(t),\qquad
 \mathcal Q\in\mathbb Z_{\ge0}[t].
\end{equation}
In particular, writing $c_q$ for the coefficient on the left and $\beta_q=[t^q]P_*(t)$, one has $c_q\ge\beta_q$ and
$\sum_{q=0}^j(-1)^{j-q}c_q\ge\sum_{q=0}^j(-1)^{j-q}\beta_q$ for every $j\ge0$.
\end{theorem}

Write $d_-,d_+$ for the fixed-point degrees in the immediately adjacent lower and upper chambers. These are the Leray--Schauder degrees of the identity minus $u\mapsto 2\kappa P_g^{-1}(\mu_u-V^{-1}dV_g)$ on the mean-zero space, with the usual identity orientation.
\begin{theorem}[Degree at resonance]\label{deg}
Under \eqref{W}, the mean-zero solution set at $\kappa_*$ is compact, and its Leray--Schauder degree is
\begin{equation}\label{degree}
 \boxed{\displaystyle d_{\fp}(\kappa_*)
   =d_-+\sum_{C(A)<0}j_A
   =d_+-\sum_{C(A)>0}j_A.}
\end{equation}
The sums extend over all representations in \eqref{wall}. If
\[
 (1-x)^{\chi(M)-N}\prod_{j=1}^N(1-x^{\mu_j})
                  =\sum_{\sigma\in\mathcal W} b_\sigma x^\sigma ,
\]
then $d_-=\sum_{\sigma<\rho_*}b_\sigma$ and $d_+=\sum_{\sigma\le\rho_*}b_\sigma$. In particular, a nonzero right-hand side in \eqref{degree} gives a solution of \eqref{normalized} at $\kappa_*$.
\end{theorem}
The series is formal and its support is locally finite; equal exponents are combined. Section~\ref{sub:degree} derives the degree from the Morse relations and allows degenerate genuine solutions.

\subsection{Relation to earlier work}
Ahmedou--Ben Ayed--Lucia \cite[\S\S4--5]{ABL} developed a
critical-point-at-infinity approach for the Dirichlet mean field
equation. Ahmedou--Ben Ayed \cite[\S2 and \S\S4--5]{ABA} treat
mixed ordinary and pinned centers for positive nonintegral conical
orders. The regular four-dimensional predecessor
\cite{regular-resonance} gives a Morse theory for the resonant
$Q$-curvature problem. For the boundary $Q$-curvature problem,
Ahmedou--Kallel--Ndiaye \cite{AKN} relate the local end reduction
to the topology of low sublevels. Ahmedou--Hu--Wang--Zhang
\cite[Proposition~1.1 and \S\S5--6]{Keller} use a complete
descending flow and relative homology at the retained ends for
the resonant Neumann problem. In the present problem the height
reduction incorporates three distinct moment regimes
from \cite[Theorem~12.6]{AWZBlowup}, and the endpoint groups
depend on the strict marked-point quota and the negative
Paneitz eigenspace.

\smallskip
\noindent\textit{Use of artificial intelligence.}
OpenAI's ChatGPT (Codex) assisted with drafting and revising proofs,
checking calculations and references, and editing the \LaTeX{} source.

\section{Variational setting and exact concentration}\label{sec:prelim}
\subsection{Exact solutions approaching the wall}
For $m<1/2$, the constant term of the regular density gives the subtraction in \eqref{finitepart}. The linear term has zero spherical average and the quadratic Taylor remainder contributes $O(r^{1-4\mu_j})\,dr$, integrable for $\mu_j<1/2$. This also permits independent deletion radii.

The exact-solution input used throughout is the following consequence of \cite[Theorems~10.2 and 12.6]{AWZBlowup}. For any sequence of solutions with $\kappa\to\kappa_*$ unbounded in $\cH$, a subsequence determines one representation in \eqref{wall}, a critical configuration $A$, and a common logarithmic height $T_0\to\infty$ such that
\begin{equation}\label{rate}
 \frac{\kappa-\kappa_*}{\kappa_*}=C(A)S_m(T_0)+o(S_m(T_0)).
\end{equation}
For the inverse peak scales $\lambda_s$ and $\tau=(\int_M K_De^{4u}\dV)^{-1}$, the matching is
$\kappa\tau h_s\lambda_s^{4\mu_s}\to1$ and $\kappa\tau e^{4T_0}\to1$ (using the moving centers in $h_s$). Only an ordinary bubble has center parameters. In the extended projections the amplitudes are measured by $(\alpha_s-1)\sqrt{\log\lambda_s}$. The matching theorem also gives $(\alpha_s-1)\log\lambda_s\to0$ on these exact sequences.

\subsection{Spectral reduction and nondegeneracy}\label{sub:spectral}
Set $s_*=\min_j(-1/\gamma_j)$, with $s_*=\infty$ in the unmarked case. The weight belongs to $L^s(M)$ for every $1<s<s_*$ and is bounded below by a positive constant. The exponential inequality on $H^2(M)$, followed by H\"older's inequality, implies that $u\mapsto d\mu_u/dV_g$ is smooth from $\cH$ to $L^s$ and maps bounded subsets into relatively compact subsets \cite[\S2]{AWZBlowup}. Indeed, after a subsequence a bounded sequence converges almost everywhere; the exponential inequality bounds every finite exponential power, and $K_D\in L^b$ for $s<b<s_*$. Jensen's inequality bounds the denominator below. Consequently solutions are in $W^{4,s}$ for every such $s$ and smooth on $X$.

For the quadratic form use the bilinear notation $Q(v,w)=\langle P_gv,w\rangle$. On the positive and negative spectral spaces use $\|w\|_+^2=Q(w)$ and $|z|_-^2=-Q(z)$, respectively. For $w\in E^+$ and $z\in E^-$, the squared Hilbert norm of $w+z$ is $\|w\|_+^2+|z|_-^2$. Direct differentiation gives
\begin{align}
 \tfrac12 DJ_\kappa(u)[v]&=Q(u,v)-2\kappa\int_M v\,d\mu_u,\nonumber\\
 \tfrac12 D^2J_\kappa(u)[v,w]&=Q(v,w)-8\kappa\Cov_{\mu_u}(v,w),\label{hessian}
\end{align}
where $\Cov_\nu(v,w)=\int vw\,d\nu-\int v\,d\nu\int w\,d\nu$ and $\operatorname{Var}_\nu(v)=\Cov_\nu(v,v)$.

\begin{lemma}\label{spectral}
For every $w\in E^+$ there is a unique maximizer $z_\kappa(w)$ of $J_\kappa(w+\cdot)$ in $E^-$. It depends $C^1$ on $w$ and is uniformly bounded for fixed $\kappa$. Critical points of $F_\kappa$ lift bijectively to solutions of \eqref{normalized}. Their nullities agree, and the indices satisfy
\begin{equation}\label{shift}
 \ind J_\kappa(w+z_\kappa(w))=k+\ind F_\kappa(w).
\end{equation}
\end{lemma}
\begin{proof}
The negative block in \eqref{hessian} is bounded above by $-2|\cdot|_-^2$ for the Hessian of $J_\kappa$. Also $|\log Z(w+z)-\log Z(w)|\le4\|z\|_\infty$. Strict concavity and coercivity from above give existence and uniqueness. At the maximum $Q(z,e)=2\kappa\int e\,d\mu_{w+z}$ for $e\in E^-$, so finite dimensionality bounds $z$. The implicit function theorem gives $C^1$ dependence, and differentiation in $w$ gives the remaining positive equation. At a critical point write the Hessian blocks in the order $E^+,E^-$ as $A,B,B^*,D$. Completing the square in the negative variable gives the direct sum of $D$ and $A-BD^{-1}B^*=D^2F_\kappa$. This proves \eqref{shift} and the nullity assertion. The estimate $0\le F_\kappa(w)-J_\kappa(w)\le C_\kappa$ follows from the same maximum formula \cite[Lemma~14.2]{AWZBlowup}.
\end{proof}

\begin{proposition}[Generic wall data and compactness]\label{compact}
At a prescribed wall, condition \eqref{W} is open and dense among smooth positive $h$ in the $C^\infty$ topology. For each such $h$ there are finitely many critical center configurations and the solution set of \eqref{normalized} at $\kappa_*$ is compact in $W^{4,s}$ for every $1<s<s_*$. On a residual subset of these coefficients every genuine solution has invertible mean-zero linearization. For these data the set of genuine solutions is finite.
\end{proposition}
\begin{proof}
There are finitely many types at fixed mass. Critical center configurations stay in a compact separated set uniformly for small $C^1$ perturbations of $\log h$, by the cluster dilation in \cite[Lemma~13.1]{AWZBlowup}. Indeed the interaction has coefficient $8$ at a colliding ordinary pair, and $4(1-\mu_j)$ or $4(1+\mu_j)$ at an inactive or active mark. Dilating the cluster differentiates the corresponding $\log(1/d)$ terms with the same strictly negative sign. Thus a Morse interaction has only finitely many critical orbits. The exclusion also gives openness by continuation of nondegenerate zeros and continuity of \eqref{allC}; for the finite part use \eqref{finitepart} and its integrable Taylor remainder.

For density vary $h$ to $he^v$. The first jets of $v$ vary the ordinary center equations independently on disjoint balls. On the compact center set finitely many such variations suffice, so finite-dimensional transversality gives Morse interactions for all types. Now keep $v$ zero to first order at every marked or critical center. The configurations remain critical under this perturbation. Their Hessians remain invertible for sufficiently small perturbations. The uniform separation bound and compactness then exclude additional critical configurations by the local inverse theorem. For $m\ge1/2$, the trace of the second jet at one dominant pole changes \eqref{allC} by a nonzero multiple of $\Delta_gv$ there. For $m<1/2$ a nonzero nonnegative bump away from all poles changes \eqref{finitepart} to first order by $\int v K_De^{4\Phi_{\varnothing,I}}\dV>0$. Finitely many nonzero linear functionals can be made simultaneously nonzero along one perturbation direction. This gives \eqref{W}, as in \cite[\S13]{AWZBlowup}.

If a sequence of genuine solutions at this fixed $h$ were unbounded in $\cH$, \eqref{rate} divided by $S_m(T_0)>0$ would give $C(A)=0$. The resulting Hilbert bound, compactness of the densities in $L^s$ and mean-zero inversion of $P_g$ give the asserted strong compactness.

It remains to arrange nondegeneracy of genuine solutions. Write $b=\log h$, and use the map to $\cH$
\[
 \mathcal S(u,b)=|P_g|^{-1}\big(P_gu-2\kappa_*(\mu_{u,b}-V^{-1}dV_g)\big),
\]
where $\mu_{u,b}$ uses $h=e^b$. Its $u$ derivative is self-adjoint Fredholm of index zero in the $|P_g|$ inner product, by \eqref{hessian}: the covariance operator is compact by finite-exponent Sobolev embedding and $d\mu_{u,b}/dV_g\in L^s$, $s>1$. If $\psi$ in its cokernel also annihilates every smooth $b$ variation $\eta$, then
$\Cov_{\mu_{u,b}}(\psi,\eta)=0$ for every such $\eta$. Hence $\psi$ is constant $\mu_{u,b}$-almost everywhere, and its volume mean zero forces $\psi=0$. Finitely many smooth variations therefore span the cokernel at each zero.

For an integer $j\ge1$, require invertibility at every zero with $\|u\|_{|P_g|}\le j$. This is open in $b$: bounded zeros for $b_\ell\to b$ have a convergent subsequence by the density compactness and elliptic inversion above, and the linearizations converge in operator norm. For density at a fixed $b$, cover its compact set of zeros in the closed $j$-ball by finitely many neighborhoods on which the joint derivative in $u$ and a common finite-dimensional space of smooth $b$ variations is onto. All zeros in that ball for sufficiently small variations belong to this union, by the same bounded compactness argument. In the union the universal zero set is a finite-dimensional manifold of dimension equal to the number of variation parameters. A regular value of its projection to those parameters is precisely one for which the $u$ derivative is onto at the corresponding zeros. Sard's theorem gives arbitrarily small such parameters. Taking the countable intersection over $j$ gives the residual assertion in the smooth topology, also after intersection with the open set \eqref{W}. Compactness at fixed wall data and isolation now imply finiteness.
\end{proof}

\section{Corrected projected bubbles and transverse variables}\label{sec:mixed}
\subsection{The extended chart and the normal variable}
Fix a representation $(p,I)$ of \eqref{wall}. On each compact subset of the configuration space we use a local ordering and take the estimates uniformly. Write $J=J_{\kappa_*}$, use $\|\cdot\|$ for the $|P_g|$ norm and $\theta_s=M_s/(2\kappa_*)$, so that $\sum_s\theta_s=1$. In exponential coordinates at $z_s$, transport the probability $\nu_s$ of \S\ref{sub:geometric} by $y\mapsto\exp_{z_s}(y/\lambda_s)$, with a fixed radial cutoff and normalization to mass one. Denote it by $\nu_s^M$. The mean-zero projected bubble satisfies
\[
 P_g\varphi_s=M_s(\nu_s^M-V^{-1}dV_g).
\]
At an active mark we use the extension of \cite[\S8]{AWZBlowup}. More explicitly, with $b=\nabla\log\Theta_s(z_s)$, let $\Psi_b$ be the entire first-harmonic correction of \cite[Proposition~4.5]{AWZBlowup}. The additional signed source is the pushforward of
\[
 \frac{M_s}{\lambda_s c_{s,\lambda_s}}\,
 \chi(y/\lambda_s)\big(b\cdot y+4\Psi_b(y)\big)\,d\nu_s(y),
 \qquad c_{s,\lambda_s}=\int\chi(y/\lambda_s)\,d\nu_s .
\]
Its mass vanishes on concentric balls by oddness. Add its mean-zero Paneitz potential to $\varphi_s$ and denote the result by $\phi_s$. At an ordinary point set $\phi_s=\varphi_s$. Differentiating an ordinary center also differentiates the $b$ of every active mark.

Let $\ell_s=\log\lambda_s$, $g_s=\mu_s\ell_s$, and
$t_s=g_s+\tfrac14\log h_s$, where $h_s$ is the matching factor of \S\ref{sub:geometric}. The corrected sum of \cite[\S8]{AWZBlowup} is $\sum_s\alpha_s\phi_s+B$ with $B\in E^-$, the amplitudes
$\zeta_s=(\alpha_s-1)\sqrt{\ell_s}$, and the actual core heights
\begin{equation}\label{actual-heights}
 H_s=t_s+2g_s(\alpha_s-1)+B(z_s).
\end{equation}
At fixed centers, $\zeta$, $B$ and $H_s$, the scales are uniquely determined for large heights by
\[
 \sqrt{\ell_s}=-\zeta_s+
 \sqrt{\zeta_s^2+(H_s-\tfrac14\log h_s-B(z_s))/\mu_s}.
\]
Thus comparability of actual heights allows different inverse length scales at different weights. The raw tangent frame uses $\phi_s/\sqrt{\ell_s}$, the logarithmic scale derivatives, the center derivatives divided by $\lambda_i$ at ordinary points, and a basis of $E^-$. Their Gram matrix is uniformly invertible by \cite[Lemma~9.5]{AWZBlowup}. Orthogonality to these full derivatives gives the local modulation coordinates.

\subsection{An affine fiber at fixed centers and scales}\label{sub:frozen}
Define the logarithmic height on $\cH$ by
\begin{equation}\label{height}
 L(u)=\frac14\log\frac{\int_M K_De^{4u}\dV}{\kappa_*},
 \qquad DL(u)[v]=\int_M v\,d\mu_u .
\end{equation}
For a prescribed large height $H$ set
\begin{equation}\label{width-variables}
 \xi_s=\frac{t_s-H}{\sqrt H},\qquad
 \alpha_s^0=1+\frac{H-t_s}{2g_s},\qquad U^0=\sum_s\alpha_s^0\phi_s .
\end{equation}
Take $|\xi|\le\epsilon$ for a fixed small $\epsilon$. Then $g_s=H+O(\sqrt H)$, $\max|\alpha_s^0-1|=O(|\xi|/\sqrt H)$, and all actual core heights of $U^0$ equal $H$. Let $E=E_{A,\ell}$ be the normal complement of the full tangent frame at \emph{unit} amplitudes and $B=0$. In particular $E\subset E^+$ and $Q(w,\phi_s)=0$ for $w\in E$. Put
\[
 f_s=\frac{\phi_s}{2g_s},\qquad
 E_\theta=\{d\in\R^n:\textstyle\sum_s\theta_s d_s=0\}.
\]
For $c\in\R$, $d\in E_\theta$, $B\in E^-$ and $w\in E$ use
\begin{equation}\label{affine}
 \begin{gathered}
 u=U^0+c\mathcal P+\mathcal D d+\mathcal V B+w,\\
 \mathcal P=\sqrt H\sum_s f_s,\qquad \mathcal D d=\sum_s d_s f_s,\qquad
 \mathcal V B=B-\sum_s B(z_s) f_s .
 \end{gathered}
\end{equation}
The common amplitude $c$ changes every core height by $c\sqrt H$; $d$ changes their relative heights and $\mathcal V B$ cancels the values of $B$ at the centers. The underlying affine fiber, with direction space $\mathrm{span}\{\phi_s\}+E^-+E$, depends only on $(A,\ell)$. In particular it is independent of the reference height $H$.

We record the uniform inverse used below. Freeze the factors normalizing the raw amplitudes, logarithmic scales and ordinary-center increments at a reference configuration. Let $p$ denote these finite parameters together with the spectral coefficients, $S(p)$ the corresponding sum of profiles and eigenfunctions, and $e_j(p)$ the normalized unit-amplitude frame, using the $|P_g|$ inner product. The equations defining the coordinates of $u$ are
\[
 {\cal R}_j(p;u)=\langle u-S(p),e_j(p)\rangle_{|P_g|}=0 .
\]
At unit amplitudes and zero remainder the parameter derivative is the negative Gram matrix. Its inverse and the normalized first derivatives of the frame are uniformly bounded by \cite[Lemma~9.5]{AWZBlowup}. At reference amplitudes $1+o(1)$ the derivative differs by $o(1)$ from this matrix. On a normalized parameter ball of radius $\delta$ the variation of $\langle\partial_iS,e_j\rangle$ is $O(\delta)$ by the same second-tangent estimates. The other derivative term is bounded by $C\|u-S(p)\|$, at most $C(\|u-S(p_0)\|+\delta)$. Preconditioning the equations by the inverse reference matrix thus gives a contraction with constant at most $1/2$, first choosing $\delta$ small and then $C\|u-S(p_0)\|<\delta/2$, so that the iteration preserves the ball. The increment and normal remainder are $O(\|u-S(p_0)\|)$, uniformly in the scales. The affine change from amplitudes to $(c,d,B)$ is invertible at fixed $(A,\ell,H)$ and respects permutations of the ordinary labels.

\paragraph{Comparison of separated representations.}
Suppose two separated representations of the same sequence have scales tending to infinity, amplitudes tending to one, bounded negative spectral components and normal remainders tending to zero in $H^2$. The atomic parts of the distributional limits of their Paneitz sources agree, which identifies the weighted supports. Indeed a bounded negative spectral component has a smooth subsequential limit, while the normal sources vanish in $H^{-2}$. Label the corresponding cores and test both sources with $\varphi(\lambda_i\exp_{z_i}^{-1}x)$ for $\varphi\in C_c^\infty(\R^4)$. These tests have bounded $H^2$ norm. The smooth terms vanish in this rescaling, and the signed source has $H^{-2}$ norm $o(1)$ by \cite[Proposition~8.1]{AWZBlowup}. The first rescaled limit is the diffuse measure $M_i\nu_i$ of total mass $M_i$.

Write $\lambda'_i$ for the competing scale. If $\lambda'_i/\lambda_i\to0$, absolute continuity of the rescaled model gives zero mass on compact sets. If the ratio tends to infinity, tightness gives a point mass or zero locally, according as the rescaled competing center has a finite limit or escapes. With comparable scales an escaping center again gives zero. A subsequential limit must therefore have a finite positive scale ratio and bounded translation. The radial model is invariant under the resulting similarity after discarding a rotation of frames. A contracting similarity, iterated, would give a point mass, and the inverse excludes expansion. A finite nonzero measure cannot be invariant under a nonzero translation. Thus the logarithmic scale differences and scaled center displacements tend to zero. By the normalized tangent bounds of \cite[Lemma~9.5]{AWZBlowup}, moving the second set of profiles to the first parameters changes their sum by $o(1)$ in $H^2$. The Gram lower bound for the amplitudes and negative modes then gives $o(1)$ differences in the natural amplitudes and spectral coefficients. Two representations in the unit-amplitude frame now lie in the same local inverse and coincide up to permutation.

\begin{lemma}[The transverse quadratic form]\label{affine-hessian}
Write $\varepsilon_H=\max_s|\alpha_s^0-1|+e^{-c_0H}$ with $c_0>0$ sufficiently small, and use the product norm $|c|+|d|+|B|_-+\|w\|$ in \eqref{affine}. At $U^0$ the half Hessian of $J_{\kappa_*}$ on this product is, up to $o(1)$ in operator norm,
\begin{equation}\label{fast-form}
 \kappa_* c^2-8\kappa_*\sum_s\theta_s d_s^2-|B|_-^2+
 \big(\|w\|^2-8\kappa_*\operatorname{Var}_{\mu_{U^0}}(w)\big).
\end{equation}
The last form is bounded below by $c_1\|w\|^2$, uniformly. On a sufficiently small fixed product ball the fiber Hessian is uniformly Lipschitz. At its origin,
\begin{align}
 |L(U^0)-H|&\le C\varepsilon_H,\nonumber\\
 |DJ_{\kappa_*}(U^0)[\mathcal P]|&\le C\sqrt H\,\varepsilon_H,\qquad
 |DJ_{\kappa_*}(U^0)[\mathcal D d+\mathcal V B+w]|
   \le C\varepsilon_H(|d|+|B|_-+\|w\|).
 \label{fast-residual}
\end{align}
\end{lemma}
\begin{proof}
On core $i$, $\phi_i-2g_i$ has bounded moments of every fixed order under the conditional model and nonlinear probabilities, and $\phi_s$ for $s\ne i$ is bounded there. Indeed the projected potential differs from $2g_i$ by $O(1+\log(1+|y|))$ on the rescaled core, and the shell masses have geometric tails \cite[\S\S8--9]{AWZBlowup}. The extended signed source has exponentially small logarithmic moments. At $U^0$ the core masses of $2\kappa_*\mu_{U^0}$ are $M_i+O(\varepsilon_H)$ and the exterior mass is exponentially small. More precisely the unnormalized integral on each core is $e^{4H}(M_i/2)(1+O(\varepsilon_H))$: apply the unnormalized test $1$ in the proof of \cite[Lemma~12.1]{AWZBlowup}. This proves the first bound in \eqref{fast-residual}.

These moment estimates are uniform under fiber perturbations. On core $i$, $\mathcal P-\sqrt H$ has conditional $L^r$ norm $O_r(H^{-1/2})$, $\mathcal D d-d_i$ has norm $O_r(|d|/H)$, and $\mathcal V B$ has norm $o(1)|B|_-$ by smoothness of $B$, the distance moment and the estimates for the $f_s$. A unit vector of $E$ has bounded exponential moments under each conditional law by \cite[Lemma~9.3]{AWZBlowup}. After removing the common constant $c\sqrt H$ from the exponent, the other amplitude variations grow on a rescaled shell of index $j$ by at most $C\epsilon(1+|j|)/\sqrt H+C\epsilon$. H\"older's inequality and the shell exponential estimate for $w$ preserve geometric tails with a smaller exponent. The exterior contributes $O(e^{-cH})$ even after multiplication by fixed powers of $H$ or by the centered exponential tilt $e^{O(\sqrt H)}$; here use the global Adams bound on $w$. These estimates prove the same conditional bounds throughout the small product ball. In particular all centered third moments of its unit directions are uniformly bounded. Differentiation of \eqref{hessian} now gives the asserted Lipschitz bound.

The kinetic Gram entries satisfy
$Q(\phi_s,\phi_s)=2g_sM_s+O(1)$ and $Q(\phi_s,\phi_i)=O(1)$ for $s\ne i$, by testing the projected sources. Pairings with a smooth negative eigenfunction are bounded. Thus $Q(\mathcal P,\mathcal P)\to\tfrac12\sum M_s=\kappa_*$, the $d$ block and its kinetic mixed entries are $o(1)$, and the $B$ block is $-|B|_-^2+o(1)|B|_-^2$. The kinetic mixed entries with $E$ vanish. In the covariance, the common constant in $\mathcal P$ drops out. The $d$ block tends to $\operatorname{Var}_\theta(d)=\sum\theta_s d_s^2$. Its mixed term with $w$ is $o(1)|d|\|w\|$: unit normal tests have zero average against each corrected model source by amplitude and negative-mode orthogonality, and their nonlinear conditional averages converge uniformly to these averages by weighted testing. The same comparison for products, and balanced coercivity \cite[Proposition~9.4]{AWZBlowup} on the unit-amplitude normal space, give the positive lower bound on the $w$ block. These estimates prove \eqref{fast-form}.

For the last two bounds in \eqref{fast-residual}, denote the core source masses by $\sigma_i$. Testing on $f_s$ gives
\[
 \tfrac12 DJ_{\kappa_*}(U^0)[f_s]
    =\alpha_s^0 M_s-\sigma_s+O(\varepsilon_H/H).
\]
Indeed the core averages of $\phi_s$ have constant terms $2g_s\mathbf1_{i=s}+c_{si}$ under the corrected source and the nonlinear law, with errors $O(\varepsilon_H)$; the bounded constants $c_{si}$ multiply $\alpha_i^0M_i-\sigma_i=O(\varepsilon_H)$, as in the amplitude calculation of \cite[Proposition~12.2]{AWZBlowup}. A smooth negative test has average $B(z_i)+O(e^{-cH})|B|_-$ on core $i$. For $w\in E$, the kinetic pairing is zero and each nonlinear core average differs from its zero corrected-source average by $O(\varepsilon_H)\|w\|$. This last estimate is the weighted source-defect estimate in the proof of \cite[Proposition~8.3]{AWZBlowup}: it uses amplitude and negative-mode orthogonality of the test, which hold on the frozen space. The marked first angular term is included in the source. Linearity proves the residual bounds.
\end{proof}

\section{Finite-dimensional reduction on height slices}\label{sec:scale}
\subsection{Stationarity at fixed height}
The normalization in \eqref{height} fixes the coefficient of $e^{-4L}$. Regularity of $L$ follows from the weighted exponential estimates of
\cite[\S2]{AWZBlowup}. Its differential can vanish only if $\mu_u=dV_g/V$.
This fixes $u$ uniquely up to a constant, and the mean-zero condition
fixes that constant. Thus all sufficiently large heights are regular.

\begin{lemma}[Stationary common-scale differential]\label{stationary}
At a regular height, $u$ is stationary for the restriction of $J_{\kappa_*}$ to $L^{-1}(L(u))$ if and only if it satisfies the exact equation
\begin{equation}\label{pde}
 P_gu=2\kappa(\mu_u-dV_g/V)
\end{equation}
for a unique real $\kappa$. On this stationary set the following equality holds as an identity of covectors on $\cH$:
\begin{equation}
 DJ_{\kappa_*}(u)=4(\kappa-\kappa_*)\,DL(u).
\end{equation}
For a sequence of these points with $\kappa\to\kappa_*$ and $L=L(u)\to\infty$ concentrating at $A$, the multiplier is
\begin{equation}\label{coef}
 4(\kappa-\kappa_*)=4\kappa_* C(A) S_m(L)+o(S_m(L)).
\end{equation}
The estimate on multipliers is uniform for points escaping in a closed isolating neighborhood of a single critical orbit of the specified type.
\end{lemma}
\begin{proof}
A covector vanishing on $\ker DL$ is a scalar multiple of $DL$. Subtract the derivatives of $J_{\kappa_*}$ and $J_\kappa$. Since $D\log Z=4DL$, the multiplier equals $4(\kappa-\kappa_*)$ precisely when $DJ_\kappa=0$. Testing against mean-zero variations gives \eqref{pde}, and both sides annihilate constants. This proves the equivalence and uniqueness. A sequence with $L\to\infty$ is unbounded in $H^2$ by Adams' inequality. With $\tau=(\int K_De^{4u}\dV)^{-1}$, the matching theorem \cite[Theorem~10.2]{AWZBlowup} gives $\kappa\tau e^{4T_0}\to1$ in its peak height $T_0$. Consequently
\[
 4(T_0-L)=\log(\kappa\tau e^{4T_0})+\log(\kappa_*/\kappa)=o(1).
\]
In each line of \eqref{order}, $S_m(T_0)/S_m(L)\to1$. the exact rate \cite[Theorem~12.6]{AWZBlowup} gives \eqref{coef}. The uniform assertion follows by taking a contrary sequence, extracting a limiting configuration using \cite[Theorem~10.2]{AWZBlowup}, and applying isolation.
\end{proof}

\subsection{Elimination at fixed logarithmic height}
We now impose $L=H$ in \eqref{width-variables}. Subscripts on height derivatives indicate whether $\ell$ or $\xi$ is fixed.

\begin{proposition}[Transverse stationary graph on a height slice]\label{height-graph}
Uniformly for $A$ in a compact separated chart, $H$ sufficiently large and $|\xi|\le\epsilon$, the equation $L(u)=H$ in \eqref{affine} determines $c$ locally in terms of $x=(d,B,w)$. On a fixed small $x$-ball the restriction of $J_{\kappa_*}$ to this slice of the affine fiber has one stationary point $u_*(A,\ell,H)$. Its parameters obey
\begin{equation}\label{slice-center}
 |d_*|+|B_*|_-+\|w_*\|+\sqrt H\,|c_{\mathrm{st}}|\le C\varepsilon_H .
\end{equation}
The constrained fiber Hessian there is invertible and has index $n-1+k$. Locally at every finite scale the energy splits, in parameter-dependent fiber coordinates, as
\begin{equation}\label{fiber-splitting}
 J_{\kappa_*}=j(A,\ell,H)-|y_-|^2+\|y_+\|^2,\qquad
 j(A,\ell,H)=J_{\kappa_*}(u_*),
\end{equation}
with $\dim y_-=n-1+k$. At fixed $H$, the critical points of $j(\cdot,\cdot,H)$ correspond exactly to constrained stationary points of the full functional in this tube.
\end{proposition}
\begin{proof}
Write $L_c=DL(u)[\mathcal P]$. On the small product ball in Lemma~\ref{affine-hessian}, $L_c=\sqrt H+O(H^{-1/2})$. The first $x$ derivative of $L$ is bounded, and the second and third derivatives involving unit product directions are uniformly bounded: these are centered moments, so the large constant in $\mathcal P$ cancels. The scalar implicit function theorem, or strict monotonicity in $c$, gives $c=\Gamma(x)$ on a fixed small ball. The mean value theorem and \eqref{fast-residual} give
\[
 |\Gamma(0)|\le C\varepsilon_H/\sqrt H,\qquad
 \|D\Gamma\|+\|D^2\Gamma\|\le C/\sqrt H .
\]
At $(c,x)=(0,0)$, \eqref{fast-residual} gives
$J_c=O(\sqrt H\,\varepsilon_H)$. The bounded Hessian gives
$J_c=O(1)$ on the product ball and changes the gradient by
$O(\varepsilon_H/\sqrt H)$ on passing to $(\Gamma(0),0)$.
Hence the slice gradient at zero is $O(\varepsilon_H)$.
On $|x|\le r$ its Hessian differs from twice the $(d,B,w)$ block of
\eqref{fast-form} by $O(r)+o(1)$. The terms containing $D\Gamma$
or $D^2\Gamma$ are $O(H^{-1/2})$, and Lemma~\ref{affine-hessian}
controls the remaining variation. The block has a uniform inverse
by its negative and positive gaps. Preconditioning by this inverse
gives a gradient iteration with Lipschitz constant at most $1/2$
on the $r$-ball, first choosing $r$ small and then $H$ large.
It preserves a ball of radius $C\varepsilon_H$ about zero. Its unique
fixed point gives \eqref{slice-center} by the bounds on $\Gamma$.

At each finite scale the local frames can be taken smooth in the ordinary centers and logarithmic scales. For the ordinary sources this follows by differentiation of the smooth translated cutoff profile. At fixed marks the pole stays fixed and scale derivatives act through $r\partial_r$ on the radial profiles; their derivatives through any fixed order obey the same local integrability bound by the radial equation and the first-harmonic formula of \cite[\S4]{AWZBlowup}. Inverting $P_g$, forming the finite Gram projection, and applying the implicit function theorem on its range therefore give smooth local coordinates and a smooth stationary graph at finite parameters. In particular the fiber changes below and their inverses have locally Lipschitz differentials there. The index follows by comparison with \eqref{fast-form}. To obtain
\eqref{fiber-splitting}, first minimize in $w$ by the uniform positive
gap. The remaining finite Hessian is the Schur complement and has a
uniform negative gap. The square roots of the positive and negative
Taylor forms give the fiber changes. On common smaller balls their
derivatives differ from the linear square-root maps by $O(r)+o(1)$:
the Taylor forms are uniformly Lipschitz in the fiber by the
Hessian bounds and the derivative of the minimizing graph. The
linear maps and their inverses are uniformly bounded by the gaps.
This gives a common radius and a uniform inverse. Varying $(A,\ell)$
with $c=\Gamma(x)$ gives a chart of $L^{-1}(H)$, so the reduced critical
equations are precisely the remaining tangent equations.
\end{proof}

\begin{corollary}[Height derivative and the limiting finite-dimensional energy]\label{slice-limit}
At any point of this graph, $\partial_H j|_{A,\ell}$ is its fiber Lagrange multiplier and is $O(\varepsilon_H)$. Put
\[
 \widehat j(A,\xi,H)=j(A,\ell(A,\xi,H),H),\qquad
 \ell_s=(H+\sqrt H\,\xi_s-\tfrac14\log h_s)/\mu_s .
\]
Uniformly for $|\xi|\le\epsilon$ on the center chart,
\begin{equation}\label{slice-energy}
 \widehat j(A,\xi,H)=c(A)+\frac12\sum_s M_s\xi_s^2+o(1),
 \qquad H\longrightarrow\infty ,
\end{equation}
with $c(A)$ given by \eqref{limitenergy}. In particular this finite-dimensional energy extends continuously to $H=\infty$ in the variables $(A,\xi,1/H)$. At a critical point of $(A,\xi)\mapsto\widehat j(A,\xi,H)$ the multiplier equals $4(\kappa-\kappa_*)$ for a solution of \eqref{pde}. Along a sequence of such points tending to infinity the limit centers are critical, $\xi\to0$, and
\begin{equation}\label{slice-slope}
 \partial_H\widehat j\big|_{A,\xi}
       =4\kappa_* C(A_*) S_m(H)+o(S_m(H))
\end{equation}
after passage to the limiting configuration $A_*$. If $C(A_*)\ne0$, then for all sufficiently large $H$ this derivative has the sign of $C(A_*)$ at every slice-critical point in a closed chart isolating $A_*$.
\end{corollary}
\begin{proof}
At fixed $(A,\ell)$ the affine fiber is the same for all $H$. Write $DJ_{\kappa_*}|_{\rm fiber}=\lambda DL|_{\rm fiber}$ at its constrained stationary point. Since $DL[\partial_Hu_*]=1$, $\partial_Hj|_{A,\ell}=\lambda$. Testing with $\mathcal P$ and using \eqref{fast-residual}, \eqref{slice-center} and the bounded fiber Hessian gives $|\lambda|\le C\varepsilon_H$. At a shape-critical point all tangent equations hold. Lemma~\ref{stationary} then gives $\lambda=4(\kappa-\kappa_*)$ and $\kappa\to\kappa_*$ as $H\to\infty$. The derivative at fixed $(A,\xi)$ equals the one at fixed $(A,\ell)$ at such a critical point.

For the energy, define $c_{ss}=\int U_s\,d\nu_s-c_s+\mu_sH_g(z_s,z_s)$ and $c_{si}=\mu_sG(z_s,z_i)$ for $i\ne s$. Pairing the corrected sources with their potentials as in \cite[proof of Proposition~12.2]{AWZBlowup} gives
\[
 Q(\phi_s,\phi_i)=M_i\big(2g_s\mathbf1_{i=s}+c_{si}\big)+o(1).
\]
The constants are uniform by separation and integrability of the logarithmic shell moments. Since $\alpha_s^0=1+O(H^{-1/2})$, all bounded terms may be evaluated at unit amplitudes. The diagonal growing term is exactly
\[
 2g_s(\alpha_s^0)^2
     =2H-\tfrac12\log h_s+(H-t_s)^2/(2g_s).
\]
Subtract $\kappa_*[4H+\log(\kappa_*/Z_0)]$ and sum. The constants combine to $c(A)$ in \eqref{limitenergy}, since
$h_s=\Theta_s(z_s)e^{-4c_s}/3$ and $\sum_sM_s=2\kappa_*$. Replacing $U^0$ by $u_*$ changes this energy by $o(1)$. Indeed $J(U^0)$ uses $L(U^0)=H+O(\varepsilon_H)$, and Taylor's formula for $J$ in \eqref{affine}, with \eqref{fast-residual} and \eqref{slice-center}, costs $O(\varepsilon_H^2)$. Finally $H/g_s\to1$ uniformly. This proves \eqref{slice-energy}.

Along a sequence of shape-critical points with $H\to\infty$ we have obtained exact solutions with $\kappa\to\kappa_*$. The matching theorem \cite[Theorem~10.2]{AWZBlowup} gives $t_s=H+o(1)$ for their modulated parameters. The comparison in \S\ref{sub:frozen} applies since $\|w_*\|=o(1)$ and the amplitudes tend to one. It identifies the two sets of logarithmic scales to $o(1)$. Since $h_s$ is continuous and positive on the separated center chart, $t_s-H=o(1)$ also on the graph, and $\xi\to0$. The exact center equation and Lemma~\ref{stationary} give \eqref{slice-slope}, uniformly in a closed isolating chart by subsequence extraction.
\end{proof}

\subsection{Entrance of exact solutions and their limiting energy}
\begin{lemma}[Finite resonant energy]\label{energy}
Let $u$ run through an escaping sequence of exact solutions with
$\kappa\to\kappa_*$ and limiting configuration $A$. The sequence
enters the height chart with $(\xi,d,B,w,\sqrt{H}\,c)\to0$, where
$H=L(u)$. Both $J_{\kappa_*}(u)$ and $J_\kappa(u)$ tend to $c(A)$.
\end{lemma}
\begin{proof}
The matching and transverse localization \cite[Theorem~10.2 and Lemma~12.5]{AWZBlowup} give $t_s=L+o(1)$ and
$\|u-\sum_s\phi_s\|=O(e^{-\delta L})$ for some $\delta>0$ at their
modulated centers and scales. The uniform inverse in
Section~\ref{sub:frozen} at this unit sum preserves these estimates:
the changes in logarithmic scales and normalized center increments,
and the new amplitudes, spectral component and normal remainder
measured from the unit sum, are $O(e^{-\delta L})$. In particular
$2g_s(\alpha_s-1)=O(\sqrt L\,e^{-\delta L})$, so
$H_s-L=\sqrt L\,c+d_s=o(1)$ by \eqref{actual-heights} and
\eqref{affine}. Its weighted mean and centered part give the asserted
entry. The exact solutions are stationary on their height slices by
Lemma~\ref{stationary}; Proposition~\ref{height-graph} and
\eqref{slice-energy} give the energy limit. The difference
$J_\kappa(u)-J_{\kappa_*}(u)$ is $o(1)$ by \eqref{coef} and
$\log Z=4L+\log(\kappa_*/Z_0)$.
\end{proof}

\subsection{Tangential differentiation on the height slice}
Write $Y=(A,\xi)$ and $f(Y,H)=\widehat j(A,\xi,H)$.
\begin{lemma}\label{lateral}
On compact separated center charts, uniformly for $|\xi|\le\epsilon$ on a smaller chart,
\begin{equation}\label{lateral-c1}
 D_Y f=D_Y f_{\partial}+o(1),\qquad
 \partial_H f\big|_{Y}=o(1),\qquad
 f_{\partial}=c(A)+\tfrac12\sum_s M_s\xi_s^2 .
\end{equation}
Consequently, in Morse coordinates $x=(x_-,x_+)$ for $f_\partial$ at $(A_*,0)$, the field $(x_-,-x_+)$ decreases $f$ at a rate bounded away from zero on each fixed closed annulus sufficiently near the origin, uniformly at large heights.
\end{lemma}
\begin{proof}
The descent and level collars use first derivatives of the stationary value. A dot denotes a variation with $|\dot A|+|\dot\xi|+|\dot H|\le1$. Put $b_0=\sqrt H\,c$ in the affine chart and keep $b_0,d$ and the coefficients of $B$ fixed for the moment. Then
\[
 \alpha_s=1+\frac{H-t_s+b_0+d_s-B(z_s)}{2g_s},\quad
 v=1+\max_s|\dot\ell_s|\le C\sqrt H,\quad
 |\dot\alpha_s|\le C H^{-1/2}.
\]
We evaluate at the stationary graph, where
$|b_0|+|d|+|B|_-+\|w\|\le C H^{-1/2}$ by \eqref{slice-center}.
On each projection cutoff support use the motion that fixes
$y=\lambda_s\exp_{z_s}^{-1}x$ in a moving orthonormal frame. Extend it by a smooth cutoff in slightly larger disjoint balls, fixing the inactive marks. At an active mark the motion is a dilation fixing the pole. Its physical velocity is $O(v)$ in $C^2$, so the derivative of the pulled-back Paneitz form is $O(v)$ in bilinear-form norm on $H^2$. Pull back the volume constraint by the same motion. To differentiate the stationary value we use a local section of the normal bundle through the stationary $w$ with
\begin{equation}\label{transport-test}
 \|\dot w\|_{H^2}\le Cv\|w\|,\qquad
 \int \dot w\,d\widetilde\nu_s=o(1)
\end{equation}
in material coordinates, where $\widetilde\nu_s$ is the corrected source of mass one before volume subtraction.

Adjoin the volume-mean row to the defining covectors of $E$ and use the $Q$ pairings in the bubble rows, which agree with the $|P_g|$ pairings on the complement of the negative modes. For the probability density $\rho_s=d\nu_s/dy$, the rescaled amplitude source is $\chi(y/\lambda_s)c_{s,\lambda_s}^{-1}\rho_s(y)\,dy$, plus the signed correction at a mark. Write $D=y\cdot\nabla_y$ and $B_b=b\cdot y+4\Psi_b$. Before pullback the scale derivative of the physical density multiplies its rescaled flat factors by $4+D$ on $\rho_s$ and by $3+D$ on $B_b\rho_s$, with the common cutoff and the derivative of its normalization. After pullback the radial factors are fixed. The signed term retains the derivative
\[
 \frac{d}{dt}\bigl(\lambda_s^{-1}B_b(y)\bigr)
   =\lambda_s^{-1}\bigl(B_{\dot b}(y)-\dot\ell_s B_b(y)\bigr)
\]
at fixed $y$, by linearity of the first-harmonic correction in the angular jet $b=\nabla\log\Theta_s(z_s)$. The same formula holds after applying $3+D$ to $B_b\rho_s$. Derivatives of the radial cutoff are supported at $|y|\asymp\lambda_s$. For an ordinary center and a center direction $e$, put $F_a(x)=\exp_a(x)$ in the chosen moving frame. The source of the $\lambda_s^{-1}$-normalized translation has rescaled density
\[
 -\operatorname{div}_y\left(
       c_{s,\lambda_s}^{-1}\chi(y/\lambda_s)\rho_s(y)
       (d_xF_a)^{-1}\partial_{a,e}F_a\big|_{x=y/\lambda_s}
                        \right).
\]
This follows by differentiating the pushforward on a test function. The vector multiplying the density is smooth in $(a,x)$ on the fixed chart, so its material derivative and that of its $y$ derivative are $O(v)$. The ordinary radial density and its first scaled spatial derivative have the same geometric shell bound. If $a=a_i$ is the translated ordinary center, the contribution to
this row from each active mark is its signed source with $b$
replaced by $\lambda_i^{-1}\partial_{a,e}b$. All these center and scale rows have mass zero, as do their material derivatives. The smooth volume and negative-mode rows have derivatives bounded by $Cv$ in $(H^2)^*$.

On shell $A_j=\{2^j<|y|<2^{j+1}\}$ the density times $|A_j|$ of each rescaled scale or translation row and its material derivative divided by $v$ is bounded by $C2^{-a|j|}$ for some $a>0$. Indeed, the corresponding bound for $\rho_s$ is $C2^{-4\mu_s|j|}$. The Euler bounds in \cite[proof of Proposition~8.1]{AWZBlowup} give the signed terms and their material derivatives divided by $v$ the additional factor $C2^j/\lambda_s$, bounded on the cutoff support. Differentiation of normalization and cutoff preserves the bound. For the amplitude source $\widetilde\nu_s$ before division by $\sqrt{\ell_s}$, the leading radial density is fixed, so its material derivative divided by $v$ has the same geometric bound times $C\min(1,2^j/\lambda_s)$, up to the derivative $O(\lambda_s^{-4\mu_s})$ of the cutoff normalization. This derivative also has mass zero. In all zero-mass pairings subtract the first-shell mean of the test and apply \cite[Lemma~9.2]{AWZBlowup}. For the extra radius factor, on the positive shells use $\min(1,r)\le r^d$ with $0<d<\min(1,a)$ and sum $(1+j)^{3/4}2^{-(a-d)j}$. The negative shells obey the same bound. This proves
\[
 \|\dot{\widetilde\nu}_s\|_{(H^2)^*}\le Cv e^{-c_1H},
 \qquad \ell_s\asymp H.
\]
The same estimates give $Cv$ for the derivatives of the centered scale and translation rows. For the normalized amplitude row the derivative of $\ell_s^{-1/2}$ contributes at most $C|\dot\ell_s|/\ell_s$ by the uniform norm of that row.

Let $T_p$ be the operator formed by these covectors on the pulled-back space, with the fixed Hilbert norm at the point of differentiation. Its Gram operator is uniformly invertible by \cite[Lemma~9.5]{AWZBlowup} and the addition of the mean row. The projection $1-T_p^*(T_pT_p^*)^{-1}T_p$ identifies nearby kernels at each finite scale. On $w\in\ker T_p$ its derivative is $-T_p^*(T_pT_p^*)^{-1}\dot T_pw$. The row bounds give the first estimate in \eqref{transport-test}. Differentiating $\int w\,d\widetilde\nu_s=0$ gives $\int\dot w\,d\widetilde\nu_s=-\int w\,d\dot{\widetilde\nu}_s=o(1)$. The shell mean of $\dot w$ is bounded since its difference from the integral against the corrected source is at most $C\|\dot w\|_{H^2}$ by the same shell estimate.

For the height constraint write $r=|y|/\lambda_s$ and divide the nonlinear density on core $s$ by $(M_s/2)\exp(4H+4b_0+4d_s)$. For $w=0$ its expression in rescaled variables is $e^{R_s}\,d\nu_s$. On a smaller fixed physical core, for some $a_0>0$ and $N<\infty$,
\begin{align}
 |R_s|&\le C\{H^{-1/2}(1+\log(2+|y|))+r+H^N e^{-a_0H}\},\notag\\
 |\dot R_s|&\le C\{H^{-1/2}(1+\log(2+|y|))+v(r+H^N e^{-a_0H})\}. \label{transport-density}
\end{align}
Here the constant part of the self-potential $2g_s$ is included in the exact height, \eqref{actual-heights}. To check the second line, write the self-potential as $2g_s+U_s(y)-c_s$ plus the regular Green term and the metric and cutoff errors, as in \cite[proof of Proposition~8.3]{AWZBlowup}. The profile $U_s(y)$ has zero material derivative. Its amplitude derivative costs $CH^{-1/2}(1+\log(2+|y|))$. The constant regular, separated and weight terms at unit amplitudes cancel against the derivative of $\log h_s$ in the extracted height. Changing their amplitudes, including the compensated $B$ term, costs $O(H^{-1/2})$ in that derivative. The remaining smooth terms have size $O(r)$ at the center, with material derivative $O(vr)$. For the logarithmic kernel this estimate before integration follows by writing the squared geodesic distance in a recentered normal chart as a smooth quadratic form in the difference of its arguments; its coefficient is the identity plus $O((r+r')^2)$ for source radius $r'$. The same bound holds for a center derivative, and a simultaneous dilation derivative costs $Cv(r+r')^2$. The Euclidean logarithmic singularity is integrated exactly by the entire profile. Integrating a positive power of $r'$ against the truncated model costs $O(H^N e^{-a_0H})$, also for the cutoff tails with a logarithmic moment. For the projected signed correction, pull both variables of its Green integral back by the same motion. On shell $j$ the signed source and $v^{-1}$ times its material derivative have density bounded by $C\lambda_s^{-1}2^j2^{-4\mu_s|j|}/|A_j|$ in rescaled coordinates, up to the cutoff $2^j\le C\lambda_s$, by \cite[Propositions~4.5 and 8.1--8.2]{AWZBlowup}. The shell average of $|G_0|$ is at most $C(1+|\log(2^j/\lambda_s)|)$ uniformly in the observation point. Its material derivative is bounded by $Cv$ because the motion distorts distances at relative rate $O(v)$ and the Green remainder is $C^1$. Summing the shells bounds the potential and $v^{-1}$ times its material derivative by $C(1+\ell_s)^2\lambda_s^{-\min(1,4\mu_s)}$. This is absorbed in the exponential errors in \eqref{transport-density}. Off the smaller cores the projected potentials and their material derivatives are bounded by $C$ and $Cv$, respectively, by separation and the fixed-annulus Green estimates. This proves \eqref{transport-density}; at a mark the singular power and the Euclidean volume dilation have already been included in $g_s$ and $d\nu_s$.

Restoring $w$ adds $4w$ to $R_s$ and $4\dot w$ to its derivative. The normalized core laws converge in local $L^{1+\eta}$ on rescaled annuli to the radial laws and have summable geometric shell bounds, for some $\eta>0$. This follows from \eqref{transport-density}, $\|w\|=O(H^{-1/2})$, and the annular exponential estimate in \cite[Lemma~9.3]{AWZBlowup}. Consequently their integrals of $\dot w$ differ by $o(1)$ from its corrected model averages: on finitely many shells use H\"older and the bounded first-shell means in \eqref{transport-test}, and on the tails sum $C(1+|j|)^N2^{-a|j|}$. For the factor $vr$, split the shell sum at $2^j=\sqrt{\lambda_s}$. Below that radius the additional factor is at most $Cv/\sqrt{\lambda_s}$. Above it the geometric tail is $O(v(1+\log\lambda_s)^N\lambda_s^{-a/2})$ after reducing $a>0$ to absorb the shell polynomial. Both bounds vanish since $v=O(\sqrt H)$ and $\log\lambda_s\asymp H$. On the exterior the contribution after division by $e^{4(H+b_0)}$ and its derivative are exponentially small, by the global exponential inequality and the bound on $\dot w$. Moving cutoff boundaries can equivalently be handled by a smooth partition in the fixed physical annuli. Thus the material derivative at fixed $b_0,d,B$ of
\[
 e^{-4(H+b_0)}\int_M K_D e^{4u}\dV
\]
is $o(1)$, while its value is bounded above and below by positive constants. Since $\partial_{b_0} L=1+O(H^{-1})$, imposing $L=H$ in this variation gives $\dot b_0=o(1)$.

Finally expand $Q(u)$ on this transported family. Its cross terms in $w$ vanish identically. The derivative of the pulled-back form and the term $2Q(w,\dot w)$ are each $O(v\|w\|^2)=o(1)$. With the constants of Corollary~\ref{slice-limit},
\[
 Q(\phi_s,\phi_i)=M_i(2g_s\one_{i=s}+c_{si})+R_{si},\qquad
 Q(\phi_s,B)=M_s B(z_s)+R_s^B,\qquad
 |R|+|\dot R|=o(1)
\]
for either remainder. Indeed pair the two transported sources against the Green kernel. In the diagonal term the flat logarithm splits into its exact scale constant and a fixed integrable logarithm of the rescaled difference. The regular terms, cutoff tails and signed terms obey the differentiated bounds just used, now integrated in both variables. In the separated term the kernel is smooth on the two cores. The same argument with the smooth function $B$ gives the second formula. The identity
\[
 2g_s\alpha_s^2
 =2(H-\tfrac14\log h_s+b_0+d_s-B(z_s))
   +\frac{(H-t_s+b_0+d_s-B(z_s))^2}{2g_s}
\]
cancels the mixed linear spectral term and, after summation, the $d$ term. The height term cancels $4\kappa_*H$ in $J$, since $\sum_sM_s=2\kappa_*$. The bounded Green and weight terms give $Dc(A)[\dot A]$ by \eqref{limitenergy}. The last square contributes $M_s\xi_s\dot\xi_s+o(1)$ after multiplication by $M_s$, since $g_s/H=1+O(H^{-1/2})$ and $|\dot g_s|\le Cv$. To identify this derivative with that of the stationary value, use the transported normal space as the fixed Hilbert fiber. At finite height the pulled-back quadratic form depends differentiably on the parameter in bilinear-form norm. The pulled-back weight has the fixed singular power at each pinned mark and a differentiable bounded multiplier there. The exponential estimates of \S\ref{sub:spectral} therefore give differentiability of the pulled-back integrals in the parameter and in the fiber variable. Eliminate $b_0$ by the scalar constraint and apply Proposition~\ref{height-graph} in this fixed fiber. Differentiation of the stationary value leaves only the partial parameter derivative just computed, because its fiber differential vanishes. This takes no spatial derivative of the normal remainder. This proves \eqref{lateral-c1}. In Morse coordinates the derivative of the limiting function along $(x_-,-x_+)$ is $-2|x|^2$, proving the final assertion.
\end{proof}

\section{Descending pseudogradient and critical points at infinity}\label{sec:flow}
Write $J=J_{\kappa_*}$. We join the local field near the selected
ends to a one-sided exterior field on the full Hilbert space.
\subsection{Monotone height on an isolated end}
Fix a critical orbit $A_*$ with $C(A_*)<0$ and a closed center chart
isolating it. The reduced energy is $f(Y,H)$.

\begin{proposition}[A local descending field]\label{local-height-field}
At all sufficiently large $H$ there is a locally Lipschitz field on
the end chart with $\dot H=1$ and
\begin{equation}\label{end-decay}
 Df[W]\le-\kappa_*|C(A_*)|S_m(H)-|\nabla_Y f|^2 .
\end{equation}
The gradient in $Y$ uses any fixed smooth local product metric. Multiplication by a strictly positive locally Lipschitz function preserves descent and the sign of $\dot H$.
\end{proposition}
\begin{proof}
Let $b(H)=\kappa_*|C(A_*)|S_m(H)$. By \eqref{slice-slope}, $f_H<-2b(H)$ at every zero of $\nabla_Y f$
for all large $H$. Set $W=\partial_H-(1+a)\nabla_Y f$, where
\[
 a=\frac{(f_H+b(H))_+}{|\nabla_Y f|^2}
 \quad\hbox{where }\nabla_Y f\ne0,\qquad
 a=0\quad\hbox{where }\nabla_Y f=0.
\]
At any point of the latter set the numerator is zero on a neighborhood. Thus $a$ is locally Lipschitz at finite heights, using $f\in C^2$ there. Direct evaluation gives \eqref{end-decay}. In the split fiber variables of \eqref{fiber-splitting}, adjoining $\dot y_-=y_-$, $\dot y_+=-y_+$ adds $-2|y_-|^2-2\|y_+\|^2$ to the energy derivative in these coordinates.
\end{proof}

\subsection{Exterior descent and entrance}
For $z=\Pi_-u$ choose a smooth function $\psi$ vanishing on $[0,M_0^2]$, with $0\le\psi'\le1$ and $\psi'=1$ for sufficiently large argument. Put
\[
 \mathcal I(u)=\kappa_*\log Z(u)+2\psi(|z|_-^2).
\]
Then $J+\mathcal I=\|u\|^2+O(1)$. In all the small separated end charts under consideration the negative projection of $u$ is bounded: pair each projected source with the finitely many smooth negative eigenfunctions and use the bounded spectral coordinates. Take $M_0$ larger than this bound.

We use the following segment construction, also on the reduced space
in Section~\ref{sec:topology}. For a segment of two continuous
gradients avoiding zero, its closest point $p$ to zero satisfies
$\langle p,q\rangle\ge\|p\|^2$ for each endpoint $q$. The constant
direction $-p/\|p\|$ pairs with both endpoints by at most $-m/2$
in a neighborhood, where $m$ is the continuous distance of the
segment to zero at the point of evaluation. A locally Lipschitz
partition of unity on the union of these neighborhoods gives a
field of norm at most one with the same inequalities.

\begin{lemma}[Descent outside the retained ends]\label{exterior56}
Assume \eqref{W} at $\kappa_*$. Fix a compact energy interval. Delete open neighborhoods of all its genuine critical points and inner tubes at large height around its configurations with $C(A)<0$ and limiting energy in the interval. For some $\delta>0$ the remaining closed part of the band supports a bounded locally Lipschitz descending field $V$ satisfying
\begin{equation}\label{controlled-height}
 D\mathcal I[V]\le\delta^{-1}(-DJ[V]),\qquad -DJ[V]>0 .
\end{equation}
The rate of descent has a positive lower bound on its intersection with any fixed norm ball. Endpoints of the energy interval are included when specifying the neighborhoods.
\end{lemma}
\begin{proof}
At $DJ+\eta D\mathcal I=0$, testing with $z$ gives
\[
 (1-2\eta\psi')|z|_-^2=-2(1-\eta)\kappa_*\int z\,d\mu_u .
\]
For $0\le\eta\le1/4$ this bounds $|z|_-$ by a constant independent of $M_0$. Enlarge $M_0$ past that constant. The equation is now exactly $DJ_{(1-\eta)\kappa_*}=0$. If there were such zeros on the closed exterior for $\eta\downarrow0$, bounded ones would converge to genuine solutions in the band. For unbounded ones the exact concentration results \cite[Theorems~10.2 and 12.6]{AWZBlowup} give a separated exact configuration with $C(A)<0$, by \eqref{rate} and \eqref{W}. The energy converges to $c(A)$ by Lemma~\ref{energy}. The unit-frame entrance proved in Lemma~\ref{energy} gives $\xi,d,B,w,\sqrt H\,c\to0$ for these exact solutions of the original-weight equation. These zeros eventually belong to the deleted inner tube. This proves exclusion for $0\le\eta\le\delta$.

Apply the segment construction above to $\nabla J+[0,\delta]\nabla\mathcal I$, writing $m(u)$ for the distance of the segment from zero. The resulting locally Lipschitz field has norm at most one and pairs with each endpoint by at most $-m(u)/2$. This gives \eqref{controlled-height}. On a bounded subset of the closed exterior $m$ is bounded below. Otherwise take $\eta_j\in[0,\delta]$ with $\nabla J(u_j)+\eta_j\nabla\mathcal I(u_j)\to0$. The nonlinear gradient terms are compact there by \S\ref{sub:spectral}, and the spectral correction is finite rank. Inverting the sign operator $1-2\Pi_-$ now gives a strongly convergent subsequence and an excluded zero.
\end{proof}

\subsection{Retained ends and the lateral faces}
Fix a retained orbit and put $c_*=c(A_*)$. Choose the Morse coordinates of Lemma~\ref{lateral} so that $f_\partial=c_*-|x_-|^2+|x_+|^2$. In the splitting \eqref{fiber-splitting} set
\[
 U=(x_-,y_-),\qquad V_+=(x_+,y_+),\qquad
 d=\dim U=4p-k_{A_*}+n-1+k .
\]
Use an isometric trivialization of the positive Hilbert fibers over the center ball and height half-line, for example by parallel transport along coordinate paths in the projection bundle. The uniform fiber radius in Proposition~\ref{height-graph} gives, for fixed small $r<R$, the closed product
\begin{equation}\label{end-block}
 {\cal N}=\{|U|\le R,\ \|V_+\|\le r,\ H\ge T\}.
\end{equation}
For sufficiently small radii and large $T$, \eqref{end-block} and a slightly larger product are injective. Otherwise take radii tending to zero and lower heights tending to infinity with competing representations. The height $L$ is intrinsic. By \eqref{slice-center} and the uniform fiber inverse, both raw normal errors tend to zero and both sets of amplitudes tend to one. The comparison in \S\ref{sub:frozen} then gives uniqueness. The same argument separates distinct orbits and types. To see that a product is closed in $\cH$, take a strongly converging sequence in it. Continuity bounds its heights. Extract convergent subsequences of its centers, scales and bounded finite-dimensional raw coefficients. Subtracting their convergent profile sums from the given sequence gives strong convergence of its normal remainders. The normal constraints pass to the limit by continuity of the projections. The split fiber map is defined on the larger tube, with a uniform margin about the smaller closed balls; equivalently its uniform Lipschitz inverse there extends to their boundaries. At the limit the continuous fiber coordinates therefore satisfy the same closed inequalities.

\begin{proposition}[Capture and exit]\label{capture56}
The exterior field can be joined to a descending field for which every product \eqref{end-block} has inward bottom and stable faces and outward unstable face. A trajectory remaining in the product has $H\to\infty$ and $(x,y_-,y_+)\to0$. Each retained orbit has such a trajectory. After leaving through its unstable face a descending trajectory cannot revisit its sufficiently small inner tube. Under \eqref{W}, every forward trajectory in a compact energy band either has precompact image or is eventually captured by a retained product. In the first case its distance to the genuine solution set tends to zero.
\end{proposition}
\begin{proof}
Fix $\epsilon_0>0$ much smaller than $R^2-r^2$. By \eqref{slice-energy} the unstable face lies below $c_*-3\epsilon_0$ for all large $T$. A concentric inner tube of radius $r_0\ll\min(r,\sqrt{\epsilon_0})$ lies above $c_*-\epsilon_0$. In the inner region use Proposition~\ref{local-height-field} with $\dot H=1$ and $\dot y_\pm=\mp y_\pm$ (thus $\dot y_-=y_-$). Use the Euclidean shape metric in the Morse coordinates. Away from this region use
\[
 \dot H=1,\qquad \dot U=U,\qquad \dot V_+=-V_+.
\]
On the annulus where $|(U,V_+)|\ge r_0$, the derivative of $J$ along this field is $-2(|U|^2+\|V_+\|^2)+o(1)$ by Lemma~\ref{lateral}. A cutoff between radii $r_0$ and $2r_0$ therefore joins descending fields. It leaves both $\dot H=1$ and the fiber equations intact. All specified face directions are strict.

Extend the local field with a cutoff equal to one on the product.
The lateral transition stays away from the inner tube, and the bottom
transition lies in a fixed finite height interval. On both regions
the local descent is uniformly strict by \eqref{end-decay} and the
annular estimate, and $D\mathcal I=4\kappa_*dH$. Thus
\eqref{controlled-height} holds, with a changed constant, on the
transitions and elsewhere outside the inner tube. Dividing the local
field by $1+\|W\|$ before patching preserves this estimate and the
face directions. On bounded height intervals the chart derivatives
and coefficients are uniformly bounded by Sections~\ref{sec:mixed}
and~\ref{sec:scale}, so the time change preserves a positive lower
bound on descent on bounded exterior sets.

Choose a norm ball whose interior contains the compact genuine
critical set $\mathcal K=\{u\in\cH:DJ(u)=0\}$.
Adams' inequality bounds $L$ on the closed ball and a slightly
larger ball. Raise the bottom heights of the end products so that
the products and their transition collars miss the larger ball. On this ball use the bounded
normalized negative gradient, with ordinary Morse fields in smaller
balls around nondegenerate genuine solutions when needed. Compactness
of bounded approximate critical sequences gives uniform strict descent
away from the genuine critical set on the ball. Since $D\mathcal I$ is
bounded there, the exterior estimate persists on the patching annulus.
Construct the exterior fields on a locally finite cover by compact
energy bands and patch them by a partition of unity in energy before
applying the fixed end and norm cutoffs. Each compact band meets only
finitely many supports, so the estimates persist. Bounded speed and
local Lipschitz regularity give a complete forward flow: at a finite
terminal time the trajectory has a limit from which it extends.

Inside \eqref{end-block} positive time changes can be undone using $H$ as time. An orbit remaining there has $H\to\infty$: on every bounded height interval its unnormalized field is bounded and the normalized height speed is bounded below. The fiber equations imply $y_-=0$ and $y_+\to0$. On any fixed shape annulus $a\le|x|\le 2a$ the shape velocity is bounded for large $H$ because $|\nabla_x f|$ is bounded below there, $f_H=o(1)$, and the coefficient in Proposition~\ref{local-height-field} is bounded. The rate $-dJ/dH$ is bounded below there, and also everywhere on $|x|\ge a$ inside the product, by \eqref{end-decay} and Lemma~\ref{lateral}. Since $J$ is bounded on the product, the orbit can neither stay outside the $a$-ball nor make infinitely many traversals of the annulus, each of which takes a fixed positive amount of height time. Let $a\downarrow0$. This proves convergence.

For existence, launch the bottom disk $\{|U|\le R,V_+=0,H=T\}$. If every trajectory exited in finite height time, its first exit would be on $|U|=R$. Transversality gives a continuous exit time, extended by zero on the boundary of the disk. Projection of the exit points onto their $U$ coordinate would retract the disk onto its sphere. For $d=0$ the exit face is empty and the single initial point is trapped. The energy gap between the exit face and the inner tube proves no return.

Consider an orbit remaining in an energy band $[a,b]$ and never
permanently captured. Each visit to an inner tube ends with an exit
from its surrounding product, after which the energy gap forbids
return to that tube. There are finitely many tubes, so the tail
avoids all of them. On each subsequent excursion outside the ball
containing $\mathcal K$, \eqref{controlled-height} bounds the increase
of $\mathcal I$ by $C(b-a)$. Its initial value is bounded on the ball
boundary by Adams' inequality, or is its value at the start of the
tail for the first excursion. Hence $J+\mathcal I=\|u\|^2+O(1)$
bounds the entire tail.

On its bounded range the patched field has uniform strict descent at
distance at least $\varepsilon/2$ from $\mathcal K$, by
Lemma~\ref{exterior56} and the annular estimates. If $\mathcal K$
is empty this excludes an orbit in the band. Otherwise the orbit
cannot remain at that distance or cross infinitely often from
distance $\varepsilon/2$ to $\varepsilon$, since bounded speed gives
each crossing a fixed energy cost. Thus its distance to $\mathcal K$
tends to zero. Compactness of $\mathcal K$ and continuity on finite
time intervals give precompactness. An orbit captured in a product
has $H\to\infty$ and is unbounded by Adams' inequality.
\end{proof}

\begin{proof}[Proof of Theorem~\ref{thm:cpi}]
Proposition~\ref{capture56} gives the complete field, its classification
and a trapped orbit for each $C(A)<0$. Along a trapped trajectory
$H\to\infty$, $(A,\xi)\to(A_*,0)$, $y_-=0$, and $y_+$ decays
exponentially in height time. The uniform fiber inverse and
\eqref{slice-center} give
$|d|+|B|_-+\|w\|=O(H^{-1/2})$.
The bounds for $\Gamma$ imply
$\sqrt H\,|c|\le C(\varepsilon_H+|d|+|B|_-+\|w\|)\to0$.
Set $v=u-w$ in \eqref{affine}. Its amplitudes tend to one and its
core heights are $H_s=H+\sqrt H\,c+d_s=H+o(1)$. With zero normal
remainder, \cite[Proposition~8.3]{AWZBlowup} gives
$\int_{B_s}K_De^{4v}\dV=e^{4H_s}(M_s/2+o(1))$.
Let $\pi_s$ be the probability obtained by normalizing this density
on $B_s$. For $w\ne0$, the unit vector $w/\|w\|$ lies in the frozen
normal space $E_{A,\ell}$. Its amplitude and negative-mode
orthogonality give zero mean against the corrected source of each
core. The shell estimate in the proof of Lemma~\ref{affine-hessian}
therefore bounds $\int e^{8|w|/\|w\|}\,d\pi_s$ uniformly. Since
$|e^t-1|\le |t|e^{|t|}$, H\"older's inequality gives
$\int|e^{4w}-1|\,d\pi_s=O(\|w\|)$, also for $w=0$.
Thus the unnormalized core mass for $u$ is $e^{4H_s}(M_s/2+o(1))$, and its
conditional law differs from $\pi_s$ by $o(1)$ in total variation.
The geometric shell tails of $\pi_s$ and $\lambda_s\to\infty$ give
concentration at the limiting center $z_s^*$. By \eqref{height} the
total mass is $\kappa_*e^{4H}$, and $\sum_sM_s=2\kappa_*$, so the
exterior fraction tends to zero. This proves
$\mu_u\rightharpoonup\sum_s M_s\delta_{z_s^*}/(2\kappa_*)$.
The energy tends to $c(A_*)$ by \eqref{slice-energy} and
\eqref{fiber-splitting}. The support identifies the marked subset
and the unordered ordinary centers.
\end{proof}

\section{Morse reduction and change of sublevel topology}\label{sec:groups}
We use coefficients in $\FF$.
\subsection{The local change of sublevels}
Let $\mathcal N$ be the product \eqref{end-block}. For a functional $g$ on it write $\mathcal N^a(g)=\mathcal N\cap\{g\le a\}$. For $a<b$ put $\mathcal N^{b,a}(g)=(\mathcal N^b(g),\mathcal N^a(g))$.
\begin{lemma}[Morse reduction at a retained end]\label{group56}
Let $\mathcal N$ be the closed end product \eqref{end-block} for a critical
point at infinity of energy $c_*$. For every sufficiently small
$\varepsilon>0$ and sufficiently large bottom height $T$, its sublevel
pair $\mathcal N^{c_*+\varepsilon,c_*-\varepsilon}(J)$ is homotopy
equivalent to $(D^d,S^{d-1})$, where $d=4p-k_{A_*}+n-1+k$.
For $d=0$ the pair is $(\{\mathrm{pt}\},\varnothing)$.
The assertion applies separately to each end at a shared limiting energy.
\end{lemma}
\begin{proof}
Write $q=c_*-|U|^2+\|V_+\|^2$ and $q_\theta=(1-\theta)J+\theta q$, $0\le \theta\le1$.
Take $\varepsilon<\min(r^2,R^2-r^2)/16$. By
\eqref{fiber-splitting}, \eqref{slice-energy} and \eqref{lateral-c1},
$J-q=o(1)$ in value and in first derivatives in $(U,V_+)$, uniformly
as $H\to\infty$. Set $\ell_\pm=c_*\pm\varepsilon$. Take $T$ large enough
that $|J-q|<\varepsilon/8$ throughout the product. The unstable
face is strictly below both level collars $|q_\theta-\ell_\pm|\le
\varepsilon/4$.

Choose a smooth cutoff $\zeta$ with $0\le\zeta\le1$, $\zeta=1$ on
$[0,r^2/2]$ and $\zeta=0$ on $[3r^2/4,\infty)$. Use the fixed-height
direction
\[
 R_0=(U,-\zeta(\|V_+\|^2)V_+,0).
\]
It is tangent to the stable and bottom faces. On either collar,
$|U|^2+\|V_+\|^2\ge\varepsilon/2$.
If $\|V_+\|^2\ge r^2/2$, the upper energy bound implies
$|U|^2\ge r^2/2-2\varepsilon$. Thus
$Dq[R_0]=-2(|U|^2+\zeta\|V_+\|^2)$ is uniformly negative
on the collars. By the tangential estimate, $Dq_\theta[R_0]$ is uniformly
negative there for all $\theta\in[0,1]$ after enlarging $T$.

Let $\eta$ be smooth, equal to one near $\ell_-$ and $\ell_+$ and
supported in their disjoint collars. Vary the chart point with $\theta$
according to
\[
 \frac{d}{d\theta}(U,V_+,H)
   =-\eta(q_\theta)\frac{\partial_\theta q_\theta}{Dq_\theta[R_0]}\,R_0 .
\]
The field is zero off the support of $\eta(q_\theta)$, tangent to the stable
and bottom faces, and zero near the unstable face. It is locally
Lipschitz at finite height and uniformly bounded in chart norm. The
height is constant, so solutions extend in both directions on $[0,1]$
by completeness of the closed product at that height and local
continuation in its larger chart. The chain rule makes $q_\theta$ constant
along the continuation near either fixed level. Uniqueness therefore
preserves the two sublevel inequalities in both directions. The
time-one map is a homeomorphism of sublevel pairs from $J$ to $q$.

For $q$, contract $V_+$ linearly to zero and slide $H$ to $T$.
These are homotopies of pairs. The remaining pair is the $U$ disk and
its closed annulus $|U|\ge\sqrt\varepsilon$, with the annulus a collar
of the sphere. When $d>0$, the radial map $U\mapsto \min(R,R|U|/\sqrt\varepsilon)\,U/|U|$, extended by zero at the origin, sends the annulus onto the sphere. Its linear interpolation with the identity preserves the annulus and fixes the sphere, giving the equivalence of pairs. If $d=0$ the annulus is empty.

In the quadratic pair, take the relative open neighborhood
$\{q<\ell_-+\delta\}$ for $0<\delta<\varepsilon/4$. Follow $R_0$
until hitting $\ell_-$ and stop, fixing the lower member. The
collar is disjoint from the unstable face and $R_0$ is tangent to
the other faces. Uniform descent gives finite, continuous hitting
times. The continuation homeomorphism transports this neighborhood
and its retraction to the pair for $J$. When the lower member is
empty, use the empty neighborhood.
\end{proof}

\subsection{Deformation and excision on the full space}
Assume \eqref{W} and nondegeneracy of the genuine solutions.
\begin{proof}[Proof of Theorem~\ref{change56}]
Fix a regular band $[a,b]$. There is no genuine critical point in the band and no retained configuration with limiting energy in it. Lemma~\ref{exterior56} therefore gives its exterior field on the entire closed band.
The controlled-height estimate bounds the increase of $\mathcal I$ by a constant
times $b-a$. The norm stays bounded in terms of the initial point because
$J+\mathcal I=\|u\|^2+O(1)$. Strict descent is uniform on the resulting bounded
part of the band. The hitting time of $J^a$ is finite and continuous
by strict crossing, and the stopped flow gives the retraction.

At $c_*$ choose disjoint products as in Proposition~\ref{capture56} for its retained ends and ordinary Morse products for its genuine critical points. The genuine products have coordinates $U,V_+$ without the height factor and $\dim U=k+i(u)$ by \eqref{shift}. Choose $l=c_*-\varepsilon$ and $h_*=c_*+\varepsilon$ so that the closed band contains no other exceptional value and each nonempty unstable face lies below $l$. Enlarge the bottom heights as in Lemma~\ref{group56}. Apply Lemma~\ref{exterior56} on this band, deleting only smaller neighborhoods inside these products. Patch with the product fields as in Proposition~\ref{capture56}, keeping them on a collar of each product boundary. Outside the products this field satisfies \eqref{controlled-height} and has uniformly positive descent on norm-bounded sets. Stop it when it first reaches
\[
 \mathcal D=J^l\ \cup\!
   \bigcup_{\mathcal N}(\mathcal N\cap J^{h_*}).
\]
Until this happens the orbit stays outside all deleted inner neighborhoods. The integrated estimate \eqref{controlled-height}, the norm bound and strict descent on bounded sets give a finite entrance time. At an entrance above $l$ the possible faces are the stable face and, at infinity, the finite bottom face. The field points strictly in on both. At their corner the time of simultaneous satisfaction of the two inequalities is locally the maximum of their individual crossing times. At $J=l$ outside the products the level is crossed strictly. At an intersection with an entrance face use the earlier of the lower-level time and the product-entrance time. The unstable face is separated in energy from this level. These facts give continuity of first entrance, including extension by zero on $\mathcal D$. The stopped flow retracts $J^{h_*}$ onto $\mathcal D$ relative to $J^l$.

Write $P_i$ for the finitely many closed upper product sublevels
and $L_i=P_i\cap J^l$. The sets $P_i$ are disjoint. For each nonempty
$L_i$, choose the relative open neighborhood of $L_i$ in $P_i$
and its stopped retraction from Lemma~\ref{group56}. In a genuine
Morse product the same construction applies directly to the
quadratic form without the height coordinate. Use the empty
neighborhood for $L_i=\varnothing$. The union of these neighborhoods
with $J^l$ is open in $\mathcal D$: its complement is the finite
union of the closed complements in the $P_i$. The retractions
glue with the identity on $J^l$ by the finite closed covering of
this neighborhood. This proves that $(\mathcal D,J^l)$ is a good
pair in the sense of \cite[Proposition~2.22]{Hatcher}. When
$J^l\ne\varnothing$, its collapse gives the finite wedge of
$P_i/L_i$ for $L_i\ne\varnothing$, together with the clopen
pieces $P_i$ for $L_i=\varnothing$. The quotient theorem and
reduced homology of the wedge give $\bigoplus_iH_q(P_i,L_i)$.
When $J^l=\varnothing$, this is the ordinary homology of the
disjoint union of the $P_i$. Applying the disk-pair computation
gives \eqref{full-local-sum}. Products are used only for $C<0$: the
exterior estimate already covers the positive-coefficient
configurations, even at the same limiting value.
\end{proof}

\section{Endpoint sublevels and strict weighted barycenters}\label{sec:topology}
We compare the resonant reduced sublevels with the adjacent lower chamber. Throughout this section we assume only $\ker P_g=\R$ and $-1<\gamma_j<0$.

For a $C^2$ functional $g$ on a Hilbert space and a regular level $c$,
the open sublevel $\{g<c\}$ includes into $\{g\le c\}$ by a homotopy
equivalence. The bounded locally Lipschitz field
$-\nabla g/(1+\|\nabla g\|)$ has a complete forward flow, preserves
both sets, and sends every boundary point into the open sublevel
at every positive time. Its time-one map and the flow homotopies
give the equivalence. We use this also for the reduced functionals,
which are $C^2$ by Lemma~\ref{spectral} and the envelope formula.

\subsection{The affine perturbation of the reduced functional}
Use the maximizer $z(w)=z_{\kappa_*}(w)$ of Lemma~\ref{spectral} and write $U(w)=Q(w)$ on $E^+$. For this fixed reduction set
\[
 K(w)=|z(w)|_-^2+\kappa_*\log Z(w+z(w)),\qquad
 G_s(w)=U(w)-sK(w).
\]
Thus $G_1=F_{\kappa_*}$, and the maximizer in $G_s$ is fixed at
$\kappa_*$. Let $\rho_{\rm prev}$ be the preceding wall weight,
with value zero below the first wall.

\begin{lemma}[Uniform critical values for the affine family]\label{aux}
Choose $0<a_0<a<1$ with $\rho_{\rm prev}<a_0\rho_*$. There is $C_0$ such that
\begin{equation}\label{mixedbound}
 DG_s(w)=0,\qquad a_0\le s\le b\le1
 \quad\Longrightarrow\quad |G_b(w)|\le C_0 .
\end{equation}
For a compact subinterval of $[a_0,1)$ the corresponding critical points are bounded in $E^+$. For every $t\in(a_0,1]$, the lower endpoint sublevel of $G_t$ is a strong deformation retract of the upper one whenever the compact interval between the endpoints is disjoint from $[-C_0,C_0]$.
\end{lemma}
\begin{proof}
Set $z=z(w)$, $u=w+z$. The two spectral equations are
$P_gw=2s\kappa_*\Pi_+(\mu_u-dV_g/V)$ and
$P_gz=2\kappa_*\Pi_-\mu_u$.
Consequently $\widehat u=w+sz$ solves \eqref{normalized} at
$\kappa=s\kappa_*$ for $\widehat K_D=K_D e^{4(1-s)z}$:
\begin{equation}\label{hat}
 P_g\widehat u=2s\kappa_*(\mu_u-dV_g/V),\qquad
 \widehat K_D e^{4\widehat u}=K_D e^{4u}.
\end{equation}
The smooth multipliers are precompact in each $C^k$ norm and uniformly positive: $z$ ranges in a bounded subset of the finite-dimensional space of smooth negative eigenfunctions. For a sequence $s\to s_\infty<1$ in the indicated interval the limiting mass is nonresonant. Nonresonant compactness \cite[Lemma~14.1]{AWZBlowup} bounds $\widehat u$, hence $w$.

For an unbounded sequence with $s\to1$, write
\[
 H_*=\tfrac14\log\left(\int_M\widehat K_De^{4\widehat u}\dV/\kappa_*\right),
 \qquad dF=2s\kappa_*\mu_u .
\]
Here $he^{4(1-s)z}\to h$ in $C^\infty$. The peak estimates
for smoothly varying coefficients in \cite[Theorem~10.2]{AWZBlowup}
apply to \eqref{hat}. Their height is exactly $H_*$, since the
two unnormalized densities in \eqref{hat} coincide. They give
$H_*\to\infty$ and
\[
 Q(\widehat u)=4s\kappa_*H_*+O(1),\qquad
 (1-s)H_*=o(1).
\]
The second bound follows also by summing the exponentially small
errors in the masses of the cores and the exterior. The use of
this energy estimate for the affine family is analogous to
\cite[\S5.1]{AN}.

Since $Q(\widehat u)=U(w)-s^2|z|_-^2$ and $\log Z(u)=4H_*+\log(\kappa_*/Z_0)$, we have $G_s(w)=O(1)$ and
\[
 G_b(w)=G_s(w)-(b-s)K(w),\qquad
 (b-s)|K(w)|\le(1-s)(4\kappa_*H_*+O(1)).
\]
This contradicts divergence of the auxiliary values and proves \eqref{mixedbound}.

For the deformation fix $t\in(a_0,1]$ and $0<\delta<t-a_0$. On a closed band of $G_t$ outside $[-C_0,C_0]$, the segment $\nabla G_t+[0,\delta]\nabla K$ excludes zero by \eqref{mixedbound}, since its points are $\nabla G_{t-\eta}$. Its distance $m(w)$ from zero is bounded below on norm-bounded parts of the band. Indeed $\nabla K$ is compact on bounded sets by the envelope formula and \S\ref{sub:spectral}, so $2w_j-(t-\eta_j)\nabla K(w_j)\to0$ for bounded $w_j$ would yield a convergent subsequence and an excluded zero.

In Lucia's one-sided deformation \cite[Proposition~1.1]{Lucia},
\eqref{mixedbound} supplies the bound at the target parameter.
The segment construction of Section~\ref{sec:flow} gives a bounded
field $W$ pairing with both endpoints by at most $-m(w)/2$. In particular $\delta DK[W]\le-DG_t[W]$. Until the lower level is reached the energy drop bounds the increase of $K$, hence also $U=G_t+tK$. Strict descent is then uniform on the bounded range of each trajectory, so the lower level is reached in finite time. Strict crossing makes the hitting time continuous, with value zero on the lower sublevel. The stopped flow gives the retraction.
\end{proof}

\subsection{The adjacent lower chamber}
The following deformation uses the energy bound of Lemma~\ref{aux}.
For the regular Paneitz and Neumann problems see \cite{AN} and
\cite[technical supplement, \S S.3]{Keller}.
\begin{proposition}\label{reinforce}
For $a$ as in Lemma~\ref{aux} and all sufficiently large $R$, the smaller sublevel in
\begin{equation}
 G_a^{-aR}=\{G_1+(a^{-1}-1)U\le-R\}
       \ \lhook\joinrel\longrightarrow\ G_1^{-R}
\end{equation}
is a strong deformation retract of the larger one.
\end{proposition}
\begin{proof}
Take $aR>C_0$ and $-R<G_1(0)$. On $G_1^{-R}$ one has $K-R\ge U>0$, and $U$ is bounded away from zero since a neighborhood of the origin misses this sublevel. Put $r=U/(K-R)$. Then $r\le1$, $G_r=-rR$, and the smaller sublevel is $\{r\le a\}$ within $G_1^{-R}$. On the strip $a\le r\le1$ set
\[
 Y=\nabla K,\qquad X=2w-rY=(K-R)\nabla r .
\]
Fix $0<\delta<a-a_0$. At each point the segment $X+[0,\delta]Y$ avoids zero. Indeed $X+\eta Y=0$ gives $DG_{r-\eta}(w)=0$ with $a_0<r-\eta\le r\le1$ and $|G_r(w)|=rR>C_0$, contrary to \eqref{mixedbound}. On every bounded part of the strip the distance $m(w)$ of this segment to zero has a positive lower bound. Otherwise $2w_j-(r_j-\eta_j)\nabla K(w_j)\to0$ for bounded $w_j$ and a subsequence of the scalar parameters converges. Compactness of $\nabla K$ on bounded sets, proved in Lemma~\ref{aux} using the envelope formula, makes $w_j$ converge strongly on a further subsequence. Its limit is in the strip and gives an excluded zero.

The segment construction of Section~\ref{sec:flow} gives a locally Lipschitz field $W$ of norm at most one pairing with each endpoint by at most $-m(w)/2$. In particular
\begin{equation}
 Dr[W]\le-\frac{m(w)}{2(K-R)},\qquad
 D\log(K-R)[W]\le-\delta^{-1}Dr[W].
\end{equation}
Integration bounds $K-R$ and $U=r(K-R)$ until $r=a$ in terms of the initial point. Thus $w$ stays bounded, and the first inequality gives strict descent at a rate bounded below there. The bounded field extends until first arrival. The hitting time is finite and continuous by strict crossing at $r=a$, with extension by zero on $r\le a$. Stopping the flow gives the claimed deformation.
\end{proof}

\begin{lemma}[Identification of the reinforced sublevel]\label{interleave}
For sufficiently large $\ell$, $G_a^{-\ell}$ has the homotopy type of a sufficiently negative sublevel of $F_{a\kappa_*}$.
\end{lemma}
\begin{proof}
Both functions differ by a uniformly bounded amount from $U-a\kappa_*\log Z(w)$. This follows from the bounded comparison in Lemma~\ref{spectral} for $F_{a\kappa_*}$ and from $G_a=(1-a)U+aF_{\kappa_*}$ for $G_a$. Fix $D>0$ larger than the bound on their difference.

All inclusions between sufficiently low sublevels of $F_{a\kappa_*}$ are homotopy equivalences by nonresonant compactness and the deformation in \cite[Lemmas~14.1 and 14.3]{AWZBlowup}. For $G_a$ the same conclusion follows from Lemma~\ref{aux} with $t=a$.

Denote $A_\ell=G_a^{-\ell}$ and $B_\ell=F_{a\kappa_*}^{-\ell}$. In the chain
\[
 A_{\ell+2D}\ \subset\ B_{\ell+D}\ \subset\ A_\ell\ \subset\ B_{\ell-D}
\]
the first two inclusions composed and the last two composed are homotopy equivalences. In the homotopy category the middle inclusion has both a left inverse and a right inverse. It is therefore a homotopy equivalence.
\end{proof}

\subsection{The weighted barycenters and the reverse homotopy}
Give an ordinary point cost $\omega(x)=1$ and a mark cost $\omega(q_j)=\mu_j$. For $\rho>0$ let
\begin{equation}
 \B_{<\rho}=\left\{\sum_{i=1}^k t_i\delta_{x_i}\,:\,
   t_i>0,\ \sum_i t_i=1,\ x_i\text{ distinct},\
                      \sum_i\omega(x_i)<\rho\right\}
\end{equation}
with the weak topology. Since the cost of a support is a wall weight, $\B_{<a\rho_*}=\B_{<\rho_*}$. We give the reverse homotopy in the nonresonant chamber.

\begin{proposition}\label{nonres}
Fix $\kappa=8\pi^2\rho>0$ with $\rho$ strictly between consecutive wall weights. If $\B_{<\rho}\ne\varnothing$, every sufficiently negative sublevel $F_\kappa^{-R}$ is homotopy equivalent to $\B_{<\rho}$. If $\B_{<\rho}=\varnothing$, these sublevels are empty.
\end{proposition}
\begin{proof}
Write $F=F_\kappa$, $B=\B_{<\rho}$. The empty case is the coercivity statement following \cite[Lemma~14.4]{AWZBlowup}. Otherwise let $b$ be the largest support cost strictly below $\rho$; then $0<b<\rho$. The number of atoms is uniformly bounded and $B$ is compact. The concentration and test-function constructions in \cite[Lemma~14.4 and Proposition~14.5]{AWZBlowup} give a retraction $\pi_B$ of a weak neighborhood onto $B$ and maps
\[
 \Psi(w)=\pi_B(\mu_w),\qquad
 \Phi_T(\sigma)=2\kappa G_+\nu_{T,\sigma},\qquad \Psi\Phi_T\simeq\operatorname{id}_B .
\]
Here $G_+$ is the inverse on $E^+$, extended by zero on constants and negative modes. The probability $\nu_{T,\sigma}$ replaces each atom at $x$ by normalized $m_D$ on the ball of $m_D$-mass $e^{-4T}$ about $x$. For any sufficiently large test depth $T$ its image lies below $-R$, and all larger depths are homotopic there. We prove the other composition.

\smallskip\emph{A margin in the energy.}
For the present nonresonant $\kappa$, write $F=U-K_\kappa$. Choose $0<t_0<t_1<1$ with $[t_0\rho,\rho]$ inside the chamber. The critical points of $U-tK_\kappa$ are uniformly bounded for $t_0\le t\le1$ by \eqref{hat}, with $\kappa$ in place of $\kappa_*$, and \cite[Lemma~14.1]{AWZBlowup}. Their values under $U-\theta K_\kappa$, $t\le\theta\le1$, are then uniformly bounded by the exponential inequality. The proof of Proposition \ref{reinforce} therefore applies here as well. For $\tau=t_1^{-1}-1>0$ and all large $D_0$ we may deform $F^{-D_0}$ in itself into
$\{F+\tau U\le-D_0\}$. Choose $R$ beyond the nonresonant critical-value bound. For every $D_0>R$, \cite[Lemma~14.3]{AWZBlowup} gives a deformation retraction of $F^{-R}$ onto $F^{-D_0}$.

For probabilities with $L^p(dV_g)$ density, $1<p<\min_j(-1/\gamma_j)$ (any fixed finite $p>1$ without marks), set
\[
 \begin{gathered}
 E(\nu)=\Ent_{m_D}(\nu)=\int_M\log(d\nu/dm_D)\,d\nu,\qquad v_\nu=2\kappa G_+\nu,\\
 \mathcal G(\nu)=\iint G_+(x,y)\,d\nu(x)d\nu(y),\\
 I_s(\nu)=(1+s)\kappa E(\nu)-4\kappa^2\mathcal G(\nu).
 \end{gathered}
\]
Compact $L^p$ families have bounded entropy since $dm_D/dV_g$ is bounded below. Jensen's entropy inequality and completion of a square give, for $y\in E^+$,
\[
 J_\kappa(y)\le I_0(\nu)+\|y-v_\nu\|_+^2,\qquad
 I_0(\mu_w)=J_\kappa(w)-\|w-v_{\mu_w}\|_+^2 .
\]
Equality in the first inequality at $(w,\mu_w)$ gives the second identity. Since $\log Z(w)$ is bounded below by Jensen and $\inf K_D>0$,
\[
 \kappa E(\mu_w)=2(w,v_{\mu_w})_+-\kappa\log Z(w)
 \le 3U(w)+\|w-v_{\mu_w}\|_+^2+C.
\]
Take $0<s_0<\min(1,\tau/3,\rho/b-1)$. On the reinforced set,
$I_{s_0}(\mu_w)\le-D_0+C$. Also the straight segment from $w$ to $v_{\mu_w}$ has $J_\kappa\le J_\kappa(w)$ by the two square identities. All these density maps are continuous in the chosen $L^p$ topology by the exponential bounds in \cite[\S2]{AWZBlowup}.

\smallskip\emph{Local contractions.}
Fix $0\le s'<s\le s_0$ and an output tolerance $\epsilon>0$ for a bounded-Lipschitz metric on probabilities. There exist $\eta>0$ and $C$ with the following property. Fix $\sigma\in B$, $D_0>0$, and a continuous family of densities with compact image in $L^p$. If every member $\nu$ satisfies $d(\nu,\sigma)<\eta$ and $I_s(\nu)\le-D_0$, the family contracts to a single density through $d(\cdot,\sigma)<\epsilon$ and $I_{s'}\le-D_0+C$.

Fix a reference barycenter $\sigma_*=\sum_i t_i\delta_{x_i}$ and choose disjoint normal balls $B_i$ about its atoms. The closures are disjoint and avoid all marks except a marked center of a ball. Choose their diameters small compared with $\epsilon$. In an ambient weak neighborhood of $\sigma_*$ their masses are close to $t_i$ and the exterior mass is at most $\alpha$, since their boundaries have zero $\sigma_*$ mass. Shrink the neighborhood for any given $\alpha>0$.

First multiply the exterior density by a scalar varying from one to zero and renormalize. This path is continuous in $L^p$ since the denominator is at least $1-\alpha$. At an intermediate probability $\nu'$ one has $\nu'\le(1-\alpha)^{-1}\nu$, $\|\nu'-\nu\|_{\rm TV}\le C\alpha$, and $E(\nu')\le(1+C\alpha)E(\nu)+C$. The last bound follows by splitting the positive and negative parts of $f\log f$, $f=d\nu/dm_D$. Since $G_+(x,y)=(8\pi^2)^{-1}\log(1/d_g(x,y))+O(1)$ and $m_D(B(x,r))\le Cr^{4\min\omega}$, entropy duality gives
\[
 \sup_x\int |G_+(x,y)|\,d\nu(y)\le C(1+E(\nu)).
\]
The same integral against $\nu'$ is bounded by domination. Polarization therefore bounds $|\mathcal G(\nu')-\mathcal G(\nu)|$ by $C\alpha(1+E(\nu))$. Choose $\alpha$ with the coefficient of $E(\nu)$ in the total error at most $\kappa(s-s')$. Thus $I_{s'}(\nu')\le I_s(\nu)+C$ throughout the restriction.

For the resulting probability $\lambda$ write $a_i=\lambda(B_i)$ and $c_i=\omega(x_i)$. Push forward in each ball by $y\mapsto e^{-t/c_i}y$ to obtain $\lambda_t$. On each bounded time interval the pushforward operators on $L^p$ are uniformly bounded by change of variables and strongly continuous by approximation with smooth compactly supported densities. Hence the path is jointly continuous in time and the input density. The relative Jacobian with respect to $m_D$ is comparable to $e^{-4t}$. Self-distances in $B_i$ scale by $e^{-t/c_i}$ up to fixed factors, and cross distances are bounded away from zero. Consequently
\[
 E(\lambda_t)\le E(\lambda)+4t+C,\qquad
 \mathcal G(\lambda_t)\ge \mathcal G(\lambda)+
                  \frac{t}{8\pi^2}\sum_i a_i^2/c_i-C.
\]
As $\sum_i a_i^2/c_i\ge1/\sum_i c_i\ge1/b$, putting $c=4\kappa(\rho/b-1-s')>0$ gives $I_{s'}(\lambda_t)\le I_{s'}(\lambda)-ct+C$.

Choose a fixed reference density supported in the same balls with masses $t_i$. At time $T$ interpolate linearly between $\lambda_T$ and the reference density dilated to time $T$. The entropy on this segment is at most $4T+C_{\mathcal F}$ by convexity, where $C_{\mathcal F}$ is finite for the given compact input family. Each intermediate mixture has total mass one, is supported in the same contracted balls, and satisfies the lower interaction bound $T/(8\pi^2 b)-C$ by the same mass inequality and the lower bound of the kernel. Thus its energy is at most $-cT+C'_{\mathcal F}$. A single finite $T$ therefore makes this entire last segment lie below $-D_0$, while the restriction and dilation estimates bound the energy on the earlier segments by $-D_0+C$ independently of the family. All three paths stay weakly close to $\sigma_*$ by the diameter and mass bounds. The endpoint is constant across the input family.

At each $\sigma_*\in B$ take an ambient weak ball of radius $2\eta_*$ on which these input estimates hold, with $\eta_*<\epsilon/4$. Choose the construction to stay within $\epsilon/2$ of $\sigma_*$. Take a finite subcover of $B$ by the balls of radii $\eta_*$, and put $\eta$ equal to the minimum of these finitely many radii. For any $\sigma\in B$ choose a covering ball containing it. Every input within $\eta$ of $\sigma$ lies in the doubled ambient ball, and every output is within $\epsilon$ of $\sigma$. Take $C$ to bound the restriction and dilation losses from the finitely many constructions. This proves the local assertion with $C,\eta$ fixed before the family and depth are chosen.

\smallskip\emph{Passage from local contractions to the reverse composition.}
Take a map from a finite complex $K$ to $F^{-R}$. Deform it to depth $-D_0$ and then to the reinforced set as above. Call the result $w$ and put $\beta=\Psi w$. As $D_0\to\infty$, the distances $d(\mu_w,\beta)$ tend uniformly to zero by concentration and continuity of $\pi_B$ and compactness of $B$. The two families to be joined over $K$ are $\mu_w$ and $\nu_{T,\beta}$. The latter is uniformly close to $\beta$ at large $T$. The test-density estimate in \cite[proof of Proposition~14.5]{AWZBlowup} gives $E(\nu_{T,\beta})\le4T$ and $\mathcal G(\nu_{T,\beta})\ge T/(8\pi^2 b)-C$. Its $I_{s_0}$ therefore tends uniformly to $-\infty$.

Set $D=\dim K+1$ and choose $s_0>s_1>\cdots>s_D=0$. Fix positive radii $\varepsilon_j$ backwards. At stage $j$, apply the local assertion with $s_{j-1},s_j$ and output radius $\varepsilon_j/2$, obtaining $\eta_j,C_j$. Choose $\varepsilon_{j-1}<\min(\eta_j/2,\varepsilon_j/2)$. These constants precede the choice of $D_0,T$. Take both depths large enough that the two endpoint families are within $\varepsilon_0$ of $\beta$ and below $-R-C_\kappa-1-\sum_jC_j$ in $I_{s_0}$. Subdivide the cylinder relative to its ends so that the oscillation of $\beta$ on each simplex is less than $\min_j\{\eta_j/2,\varepsilon_j/2\}$. At free vertices use the test family.

Suppose the boundary of a $j$-simplex has been defined. Every density there lies within $\varepsilon_{j-1}$ of $\beta$ at the corresponding point and hence within $\eta_j$ of $\beta$ at a chosen vertex. Since $E\ge0$, faces of smaller dimension and prescribed endpoint faces satisfy the required bound in $I_{s_{j-1}}$ as well. The boundary family is compact in $L^p$. Apply its contraction and descend it to the cone on the boundary. The constant endpoint makes the cone extension continuous at the apex, and the initial map fixes all boundary values. The extension has distance less than $\varepsilon_j$ from $\beta$ on the simplex and increases the upper energy bound by at most $C_j$. Finite closed gluing gives a continuous $L^p$ cylinder with $I_0\le-R-C_\kappa$.

Applying $\nu\mapsto v_\nu$ sends the cylinder into $F^{-R}$ by the square inequality and $F\le J_\kappa+C_\kappa$. Its first end joins to $w$ by the straight segment already estimated; choose $D_0>R+C_\kappa$. Its second end is $\Phi_T\beta$. Varying a sufficiently large test depth and the initial deformation changes neither homotopy class in $F^{-R}$. In particular, fix a single admissible test depth $T_*$ for the level $-R$. The cylinder constructed for a given finite complex may require a larger $T$. The homotopy $\Phi_t$, $T_*\le t\le T$, takes values in $F^{-R}$, so the original map is freely homotopic there to its composition with $\Phi_{T_*}\Psi$. For a point this gives a bijection on components together with $\Psi\Phi_{T_*}\simeq\operatorname{id}_B$. On each component take a basepoint $x=\Phi_{T_*}(b)$. If a based sphere maps trivially under $\Psi$, the cylinder makes it freely null-homotopic, hence based null-homotopic by change of basepoint along the null homotopy. Conversely $(\Psi\Phi_{T_*})_*:\pi_j(B,b)\to\pi_j(B,\Psi x)$ is an isomorphism by the homotopy on $B$. It factors through $\Psi_*:\pi_j(F^{-R},x)\to\pi_j(B,\Psi x)$, which is therefore onto for every $j\ge1$. The regular-level deformation at the start of this section gives the homotopy type of the open sublevel. The interior is covered by Hilbert balls whose nonempty finite intersections are convex, so it has CW type by \cite[Corollary~4G.3]{Hatcher}. The compact $B$ is triangulable by the moment embedding in \cite[\S14]{AWZBlowup}. Whitehead's theorem on each component now gives the equivalence \cite[Theorem~4.5]{Hatcher}.
\end{proof}

\subsection{The resonant sublevel and the negative spectral sphere}
\begin{proof}[Proof of Theorem~\ref{main}]
Proposition \ref{reinforce}, Lemma \ref{interleave} and Proposition \ref{nonres} give (i), because $a\rho_*$ and $\rho_*$ have exactly the same strictly admissible supports. If $B=\varnothing$, lower-chamber coercivity and the bounded comparison in Lemma~\ref{interleave} make $G_a^{-aR}$ empty for large $R$. Proposition~\ref{reinforce} then forces $G_1^{-R}$ to be empty as well.

For (ii), at fixed $w$ use the maximizing $z(w)$ and the negative norm on $E^-$. In the ray $z(w)+r\xi$, $|\xi|_-=1$, use the new radius
\[
 t=\bigl[F_{\kappa_*}(w)-J_{\kappa_*}(w+z(w)+r\xi)\bigr]^{1/2}.
\]
The square has value and first derivative zero at $r=0$. Its second derivative is
$2+16\kappa_*\operatorname{Var}_{\mu_{w+z(w)+r\xi}}(\xi)$ and lies between $2$ and a fixed constant by finite dimensionality. Thus $r\mapsto t$ is strictly increasing onto $[0,\infty)$, uniformly comparable to $r$, with a continuous parameter-dependent inverse. This is a fiber homeomorphism identifying $J_{\kappa_*}$ with $F_{\kappa_*}(w)-|y|^2$ on $E^+\times\R^k$.

Lemma \ref{aux} bounds all genuine critical values for $G_1=F_{\kappa_*}$, and full critical points are its lifts. For $c=-R$ sufficiently negative the two open sublevels therefore include into the closed ones by homotopy equivalences. In the open full sublevel use the cover
\[
 \{F_{\kappa_*}(w)<c\},\qquad
 \{y\ne0,\ F_{\kappa_*}(w)-|y|^2<c\}.
\]
Writing $A=\{F_{\kappa_*}<c\}$, the first set is $A\times\R^k$. The second retracts along rays to the continuous radius $1+\sqrt{(F_{\kappa_*}(w)-c)_+}$ and has type $E^+\times S^{k-1}$. The intersection is $A\times(\R^k\setminus0)$. After the evident equivalences its two maps are the projections $A\leftarrow A\times S^{k-1}\to S^{k-1}$, since $E^+$ is contractible. By \cite[Propositions~4G.1--4G.2]{Hatcher} the paracompact open union has the homotopy type of this homotopy pushout. The double mapping cylinder of the two projections is $A*S^{k-1}$. If $A=\varnothing$ only the second set occurs. This proves (ii).

Stability of the reduced inclusions follows from Lemma~\ref{aux} with $t=1$. In the full space the radial coordinate just constructed is independent of $c$. Under a level inclusion the first open cover pieces and their intersections induce the reduced level inclusion on the base, while on the second pieces projection to $(w,y/|y|)$ is a homotopy equivalence at either level. Invariance of the homotopy pushout gives the full assertion for open sublevels, and the open-to-closed equivalences give it for closed ones by the commuting inclusion square.
\end{proof}

\begin{corollary}[The lower term in the Morse relations]\label{lower-term}
Over any coefficient field $\mathbb F$, for $R\ge R_0$ and $B\ne\varnothing$,
\[
 H_q(\cH,J_{\kappa_*}^{-R};\mathbb F)
       \cong\widetilde H_{q-k-1}(B;\mathbb F).
\]
If $B=\varnothing$, the only nonzero relative group is $\mathbb F$ in degree $k$. For $\rho_*\le1$, $B$ is the realization of the simplicial complex with nonempty faces $I$ satisfying $\sum_{i\in I}\mu_i<\rho_*$.
\end{corollary}
\begin{proof}
The ambient space is contractible. For $k>0$ and $B\ne\varnothing$, the join with $S^{k-1}$ is the $k$-fold suspension. The long exact sequence of the pair gives the asserted shift, also for $k=0$. At the first wall use the sphere in Theorem \ref{main}(ii), or the empty sublevel for $k=0$. If $\rho_*\le1$ every admissible support consists of marks only, and its barycentric weights parametrize precisely the stated simplex.
\end{proof}

\subsection{The upper sublevel and the endpoint polynomial}

\begin{proposition}\label{high}
For every sufficiently large $R$, the inclusions $F_{\kappa_*}^R\hookrightarrow E^+$ and $J_{\kappa_*}^R\hookrightarrow\cH$ are homotopy equivalences. For such $R$ the full endpoint pair has Poincar\'e polynomial $t^k P_*(t)$, with $P_*$ as in \eqref{pol}.
\end{proposition}
\begin{proof}
Choose $R>C_0$ from Lemma~\ref{aux} with $0\in F_{\kappa_*}^R$. For each $b>R$, Lemma~\ref{aux} with $t=1$ retracts $F_{\kappa_*}^b$ onto $F_{\kappa_*}^R$. A path in $E^+$ joining two points of $F_{\kappa_*}^R$, or a disk filling a sphere in this sublevel, has compact image in some $F_{\kappa_*}^b$. The retraction preserves the endpoints or boundary. Thus $F_{\kappa_*}^R$ is path connected and weakly contractible. It has CW type by the regular-level and ball-cover argument in Proposition~\ref{nonres}. Whitehead's theorem gives its contractibility.

For $k>0$ apply the open cover in the proof of Theorem~\ref{main}(ii) at $c=R$. The open reduced sublevel has the same type as $F_{\kappa_*}^R$ by the regular-level deformation. It is nonempty and contractible, so its join with $S^{k-1}$ is contractible. The full regular-level deformation gives the assertion for $J_{\kappa_*}^R$. For $k=0$ the two functionals agree. The long exact sequences of pairs give an isomorphism from $H_*(J_{\kappa_*}^R,J_{\kappa_*}^{-R})$ to $H_*(\cH,J_{\kappa_*}^{-R})$. Apply Corollary~\ref{lower-term}. This endpoint calculation applies at every wall, independently of the nondegeneracy assumptions used to compute the local groups in Theorem~\ref{change56}.
\end{proof}

\begin{proposition}[Endpoint ranks for walls of weight at most one]
Suppose $N>0$, $0<\rho_*\le1$, and label a mark of smallest weight by $1$. Then
\begin{equation}\label{quota-pol}
 P_*(t)=
 \sum_{\substack{I\subset\{2,\ldots,N\}\\
                  \rho_*-\mu_1\le w(I)<\rho_*}}t^{|I|}.
\end{equation}
For $B\ne\varnothing$, the strict barycenter space is homotopy equivalent to a wedge with one sphere $S^{|I|-1}$ for each set in the sum. An empty wedge is a point.
\end{proposition}
\begin{proof}
At the first wall $\rho_*=\mu_1$, only $I=\varnothing$ occurs and \eqref{pol} gives $1$. Otherwise $B$ is the scalar quota complex of Corollary~\ref{lower-term}, with minimal vertex $1$. By \cite[Theorem~1.3]{PW} its sphere summands correspond to nonempty faces $I$ omitting $1$ with $\rho_*-\mu_1\le w(I)<\rho_*$. Their dimensions are $|I|-1$, including dimension zero. The relative shift in \eqref{pol} gives \eqref{quota-pol}.
\end{proof}

For $\rho_*>1$ the ordinary centers also enter the endpoint groups. Here is a finite chain complex for their computation. Choose a finite geometric triangulation of $M$ containing the marks as vertices, with realization in $\R^l$. Let $d$ bound the number of admissible atoms, for example $d=\lfloor\rho_*/\min_M\omega\rfloor$. Map $\sigma$ to the vector of its integrals of coordinate monomials of total degree at most $2d-1$. This vector determines $\sigma$: on the union of two supports one can separate any one point from all the others by a product of at most $2d-1$ affine functions. Thus the moment map embeds the compact $B$. Its image is semialgebraic, since each admissible support type is parametrized by distinct points in simplices and positive coefficients of sum one, with a polynomial moment map \cite[\S14]{AWZBlowup}. Take a finite triangulation of the image compatible with these finitely many strata. If $n_i$ is its number of $i$-simplices and $\partial_i$ its augmented simplicial boundary over $\FF$, then $\beta_0=0$ and
\[
 \beta_q=n_{q-1}-\operatorname{rank}\partial_{q-1}
                     -\operatorname{rank}\partial_q,\qquad q\ge1 .
\]
Here $\partial_0$ is the augmentation. Using the triangulation of the entire image incorporates faces where ordinary centers collide, reach a mark or lose a coefficient. For $N=0$ and $\rho_*=m\in\mathbb N$, $B$ is the ordinary space of barycenters with at most $m-1$ atoms. At $m=2$, $B=M$; in particular for $M=S^4$ one obtains $P_*(t)=t^5$.

\begin{lemma}[Euler characteristic of the lower chamber]\label{endpoint-euler}
With the chamber coefficients defined in Theorem~\ref{deg}, one has
\begin{equation}\label{lower-euler78}
 P_*(-1)=1-\chi(B)
   =\sum_{p+w(I)<\rho_*}(-1)^{p+|I|}
                      \binom{\chi(M)-N}{p}
   =\sum_{\sigma<\rho_*} b_\sigma=d_- .
\end{equation}
Here $\chi(\varnothing)=0$, and the generalized binomial coefficient at $p=0$ is one.
\end{lemma}
\begin{proof}
The first equality follows from \eqref{pol}. If $B=\varnothing$,
only $b_0=1$ enters the partial sum. Otherwise choose a quota
$q$ strictly between the preceding wall weight and $\rho_*$.
The non-strict quota $q$ gives exactly $B$.
Kallel's generating-series formula
\cite[Corollary~1.8]{Kallel} computes $1-\chi_c(B)$ as
$\sum_{\sigma\le q}b_\sigma$. Since $B$ is compact, $\chi_c(B)=\chi(B)$.
The exponents up to $q$ are exactly those below $\rho_*$.
Expanding the generating series gives the indicated binomial sum,
and the chamber degree formula \cite[Theorem~1.1]{AWZBlowup} gives the
last equality.
\end{proof}

\section{Morse inequalities and the resonant degree}\label{sec:graded}
Assume \eqref{W}. For $q\ge0$ put
\[
 a_q=\sum_{p+w(I)=\rho_*}
  \#\{A\in\Crit\F_{p,I}: C(A)<0,\ 4p-k_A+p+|I|-1=q\},
 \qquad \beta_q=[t^q]P_*(t).
\]
Write $A_-(t)=\sum_{q\ge0}a_qt^q$. When all genuine solutions are nondegenerate, let $m_q$ count those of reduced index $q$. Thus their full index is $k+q$. All these sequences have finite support by Proposition~\ref{compact}, \eqref{pol} and finite triangulability of $B$.

\subsection{The filtration and the strong inequalities}
\begin{proof}[Proof of Theorem~\ref{graded-morse}]
Fix $R$ as in Proposition~\ref{high}, above the absolute values of the exceptional levels in Theorem~\ref{change56}. Between $-R$ and $R$ choose regular values defining a filtration $X_0=J^{-R}\subset\cdots\subset X_s=J^R$ with one exceptional value per step, omitting regular bands by deformation. If there are no exceptional values, the endpoint groups vanish and the relation follows from the same regular deformation. With $\FF$ coefficients write
\[
 L_{j,h}=H_h(X_j,X_{j-1}),\qquad V_{j,h}=H_h(X_j,X_0),\qquad
 r_{j,h}=\operatorname{rank}\bigl(L_{j,h}\xrightarrow{\partial}V_{j-1,h-1}\bigr).
\]
The exact sequence of the triple gives
\[
 \dim L_{j,h}+\dim V_{j-1,h}-\dim V_{j,h}
               =r_{j,h}+r_{j,h+1}.
\]
This argument uses the full direct sum \eqref{full-local-sum} at a shared critical value. The groups $L_{j,h}$ vanish for $h<k$, hence so do $V_{j,h}$ by induction, and $r_{j,k}=0$. In degree $k+q$ the sum of the local ranks is $m_q+a_q$ by Theorem~\ref{change56}. The terminal rank is $\beta_q$ by Proposition~\ref{high}. Summation over $j$ gives
\begin{equation}\label{rank-identity78}
 m_q+a_q-\beta_q=e_q+e_{q-1},\qquad
 e_q=\sum_{j=1}^{s}r_{j,k+q+1}\ge0,\qquad e_{-1}=0 .
\end{equation}
Multiplication by $t^q$ and summation give \eqref{morse}, with $\mathcal Q(t)=\sum_q e_qt^q$. The coefficient inequalities follow from \eqref{rank-identity78}, and the alternating sums through $q$ have excess precisely $e_q$. This proves the theorem.
\end{proof}

\begin{corollary}[Bounds at prescribed indices]\label{index-bounds78}
For nondegenerate genuine solutions, write $x_+=\max(x,0)$. Then
\begin{align}
 m_q&\ge(\beta_q-a_q)_+,\label{weak78}\\
 \sum_{q=0}^j(-1)^{j-q}m_q
   &\ge \sum_{q=0}^j(-1)^{j-q}(\beta_q-a_q),\qquad j\ge0. \label{strong78}
\end{align}
For $\rho_*\le1$ and $N>0$, $\beta_q$ in these inequalities is the number of $q$-element subsets in \eqref{quota-pol}.
\end{corollary}
\begin{proof}
Use $m_q\ge0$ and \eqref{rank-identity78} for the first bound. In the second the difference between the two sides telescopes to $e_j$.
\end{proof}

\subsection{The Euler identity}
Put $D_*=P_*(-1)-A_-(-1)$. Evaluation of \eqref{morse} at $t=-1$ and Lemma~\ref{endpoint-euler} give
\begin{equation}\label{graded-euler78}
 \sum_q(-1)^q m_q
 =D_*=1-\chi(B)-\sum_{C(A)<0}(-1)^{i_\infty(A)}
 =d_-+\sum_{C(A)<0}(-1)^{n+k_A}.
\end{equation}
Indeed $i_\infty=4p-k_A+n-1$ gives $(-1)^{i_\infty}=-(-1)^{n+k_A}$.

\subsection{Degree from the Morse inequalities}\label{sub:degree}
On $\cH$ put $\mathcal A=|P_g|^{-1}P_g=1-2\Pi_-$ and
$\mathcal L=(\kappa_*/2)\log Z$. In the $|P_g|$ inner product,
\[
 \nabla(J/2)=\mathcal A u-K_0(u),\qquad K_0=\nabla\mathcal L,\qquad
 \mathcal T_{\kappa_*}=\mathcal A\nabla(J/2).
\]
Here $K_0$ is weak-to-strong continuous on bounded sequences by the
density compactness in Section~\ref{sub:spectral}. The fixed-point
degree on $\cH$ agrees with its degree on the mean-zero $W^{4,s}$ space,
$1<s<s_*$, by commutativity for the compact factor through this space
\cite[\S11]{AWZBlowup}. At a nondegenerate solution the local gradient sign
is $(-1)^{k+i(u)}$. Multiplication by $\mathcal A$ contributes
$(-1)^k$, giving the fixed-point sign $(-1)^{i(u)}$.

\begin{lemma}[Bounded Morse perturbation]\label{wall-morse-perturbation}
Assume \eqref{W}. There is a smooth perturbation $\widetilde J$ of $J$
on a bounded norm ball with finitely many and only nondegenerate
critical points. The sufficiently deep and high sublevels and all
the local end pairs agree with those for $J$. The signed count of
the new critical points in reduced indices, that is, with sign
$(-1)^{\ind\widetilde J-k}$, equals $d_{\fp}(\kappa_*)$.
\end{lemma}
\begin{proof}
Choose $R_0<R_1$ so that the compact critical set of $J$ lies in
$\|u\|<R_0$. On the closed annulus $R_0\le\|u\|\le R_1$,
$\|\nabla J\|$ has a positive lower bound by bounded
approximate-critical compactness. Let $\chi(\|u\|^2/2)$ equal one
on the inner ball and zero near and outside the outer sphere.
Let $P_N$ be increasing finite-rank orthogonal spectral projections
including $E^-$ and converging strongly to the identity. Set
\[
 H_{N,v}(u)=\mathcal L(u)-\mathcal L(P_Nu)
                 -\langle v,P_Nu\rangle_{|P_g|},
 \qquad
 \widetilde J/2=J/2+\chi H_{N,v},\qquad v\in P_N\cH .
\]
On the bounded outer ball $\mathcal L-\mathcal L\circ P_N\to0$ uniformly in value and first derivative.
Indeed $K_0(u)-K_0(P_Nu)\to0$ uniformly there: otherwise choose
$u_N\rightharpoonup u$, so $P_Nu_N\rightharpoonup u$ and both
images converge strongly to $K_0(u)$. The finite-rank tails
$(1-P_N)K_0(P_Nu)$ go uniformly to zero by compactness of the
image on bounded sets. The value estimate follows by integrating
the two gradients on the segment from $0$ to $u$.
Choose $N$ large and then $\|v\|$ small so the homotopy
$J/2+t\chi H_{N,v}$, $0\le t\le1$, has no zeros of its gradient
on the closed annulus. Its gradient is
\begin{align*}
 \nabla(J/2+t\chi H_{N,v})
   &= (\mathcal A+b_t(u)I)u
       -(1-t\chi)K_0(u)-t\chi(P_NK_0(P_Nu)+v),\\
 b_t&=t\chi'H_{N,v}.
\end{align*}
In the inner ball $\chi=1$. At a zero, projection onto the orthogonal
complement of $P_N\cH$ gives
$(1-P_N)u=(1-t)(1-P_N)K_0(u)$, since $\mathcal A$ is the identity
there. At $t=1$ this forces $u\in P_N\cH$, where the remaining
equation is $\mathcal A u-P_NK_0(u)=v$. Sard's theorem permits an
arbitrarily small regular value $v$ of this finite-dimensional map.
At each such point the Hessian of $\widetilde J/2$ is the direct sum
of the invertible finite block and the identity on the complement.
The critical points form a closed discrete subset of the compact
inner ball in $P_N\cH$, hence are finite. At each zero $u=P_Nu$, the Hessian of $\widetilde J/2$ restricted to $E^-\subset P_N\cH$ agrees with that of $J/2$ at $u$ and is negative definite by \eqref{hessian}. Thus each full index is at least $k$.

The same smallness choice gives $|b_t|<1/2$. Multiplication by
$(\mathcal A+b_t I)^{-1}$ makes the displayed gradient a homotopy
of identity minus compact maps on bounded sets. The inverse
multiplies the two spectral components by bounded scalars and thus
preserves relative compactness of the nonlinear terms.
At $t=0$ the map is $\mathcal T_{\kappa_*}$. At each endpoint zero
$b_1=0$ locally and the orientation computation above gives sign
$(-1)^{\ind\widetilde J-k}$. Homotopy invariance on the outer ball
and excision give the signed count.

The exact expression for the normalized gradient also gives bounded
approximate-critical compactness for $\widetilde J$. Outside the ball
$\widetilde J=J$, so the exterior estimates and the local end fields
continue to apply there. Choose the bottom heights outside a larger
norm ball. On the remaining bounded region use a normalized negative
gradient, with ordinary Morse fields in the critical products.
Away from their neighborhoods the descent rate has a positive lower
bound by bounded approximate-critical compactness. Since $D\mathcal I$ is
bounded there, the controlled-height estimate used in
Theorem~\ref{change56} persists on that region after patching.
Its stopped deformation and relative homology proof thus apply
to $\widetilde J$ with the same end groups. Finally $J$ and
$\widetilde J$ are uniformly bounded on the perturbation ball.
They have identical sublevels at sufficiently large positive and
negative levels. This proves the lemma.
\end{proof}

\begin{proof}[Proof of Theorem~\ref{deg}]
If $J$ is Morse, the signed fixed-point count is the left side of
\eqref{graded-euler78}. Otherwise the same filtration for
$\widetilde J$ in Lemma~\ref{wall-morse-perturbation} has the
original end and endpoint groups. Evaluation at $t=-1$ gives
the Euler identity with the signed count of the perturbed points,
which equals the original degree by the lemma. Hence
$d_{\fp}(\kappa_*)=d_-+\sum_{C(A)<0}j_A$ in either case.

By the configuration-space index sum \cite[Proposition~13.2]{AWZBlowup},
\[
 \sum_{(p,I),A}j_A
   =\sum_{p+w(I)=\rho_*}(-1)^{p+|I|}\binom{\chi(M)-N}{p}
   =b_{\rho_*}.
\]
Indeed the factor $p!$ from the ordered ordinary centers cancels
against passage to unordered configurations. Since
$d_+=d_-+b_{\rho_*}$ by the chamber formula, this also gives the
upper-chamber expression in \eqref{degree}. The existence assertion
follows from the existence property of the Leray--Schauder degree.
\end{proof}

\subsection{Multiplicity and existence}
\begin{corollary}[Multiplicity and existence]\label{count78}
Let $L_{\mathrm e}$ and $L_{\mathrm o}$ be the sums of $(\beta_q-a_q)_+$ over even and odd $q$, respectively. With nondegenerate genuine solutions their number $N_{\rm sol}$ satisfies
\begin{equation}\label{parity78}
 N_{\rm sol}\ge \max\{2L_{\mathrm e}-D_*,\,2L_{\mathrm o}+D_*\},
 \qquad N_{\rm sol}\equiv D_*\pmod 2 .
\end{equation}
Under \eqref{W} alone, a genuine solution exists whenever
\begin{equation}\label{obstruction78}
 A_-(t)-P_*(t)\ \notin\ (1+t)\mathbb Z_{\ge0}[t].
\end{equation}
\end{corollary}
\begin{proof}
If $E$ and $O$ are the even- and odd-index counts, then $E-O=D_*$ by \eqref{graded-euler78}. By \eqref{weak78}, $E\ge L_{\mathrm e}$ and $O\ge L_{\mathrm o}$. Use $E+O=2E-D_*=2O+D_*$ to get \eqref{parity78}. For the last assertion suppose the solution set is empty. The nondegeneracy hypothesis of Theorem~\ref{graded-morse} then holds vacuously and \eqref{morse} gives $A_--P_*=(1+t)\mathcal Q$, contradicting \eqref{obstruction78}.
\end{proof}

In particular, when $D_*=0$, any index with $\beta_q>a_q$ gives existence under \eqref{W} alone and at least two solutions for nondegenerate genuine solutions. There is also a finite test when $D_*=0$. Then $A_--P_*$ is divisible by $1+t$ in $\mathbb Z[t]$. If the difference is nonzero and has degree $h$, the coefficient of $t^j$ in the quotient is $\sum_{q=0}^j(-1)^{j-q}(a_q-\beta_q)$ for $0\le j<h$, by coefficient comparison. Any negative coefficient gives \eqref{obstruction78}, even when $a_q\ge\beta_q$ for all $q$.

\subsection{The first marked wall and a coincident wall}
\emph{The first marked wall.} Suppose $N>0$ and $\rho_*=\min_j\mu_j$. Let $s_-$ count the marks of this weight with $C<0$. Then $P_*=1$ and $A_-=s_-$, since each retained end has $p=0$, $n=1$. In degrees zero and one, \eqref{weak78} and \eqref{strong78} give
\[
 m_0\ge(1-s_-)_+,\qquad m_1\ge m_0+s_--1,\qquad D_*=1-s_- .
\]
In particular $s_-=0$ forces a solution of reduced index zero, while $s_-\ge2$ forces at least $s_--1$ distinct solutions of reduced index one, assuming genuine nondegeneracy.

\emph{Two marks of weight one half.} Suppose there are exactly two marks, $\mu_1=\mu_2=\tfrac12$, and $\rho_*=1$. Formula~\eqref{quota-pol} gives $P_*=t$. There are two wall representations. Let $\eta\in\{0,1\}$ record whether the marked pair has $C<0$, and let $o_h$ count the retained critical points of $\F_{1,\varnothing}$ of Morse index $h$. The marked pair has reduced index one and the ordinary configurations have reduced index $4-h$. Thus
\[
 A_-(t)=\eta t+\sum_{h=0}^4 o_h t^{4-h},\qquad
 D_*=-1+\eta-\sum_{h=0}^4(-1)^h o_h .
\]
The degree-one strong inequality reads
\[
 m_1\ge m_0+1+o_4-\eta-o_3 .
\]
For example, if $\eta=o_3=0$ it forces a solution of reduced index one for Morse data. If in addition $D_*=0$, there is also an even-index solution by \eqref{graded-euler78}. The two rates in \eqref{order} are $e^{-2H}$ for the ordinary type and $He^{-4H}$ for the marked pair. Their limiting energies may coincide. Independently, the chamber series is $(1-x)^{\chi(M)-2}(1-x^{1/2})^2$. Its coefficients of $x^{1/2}$ and $x$ are $-2$ and $3-\chi(M)$, respectively. Thus $d_-=-1$ and $d_+=2-\chi(M)$, which fixes the reference degrees on both sides of this wall.

\section{The resonant degree and existence for the geometric equation}\label{sec:geometric-applications}

\subsection{The two geometric masses}
Write $v_D^0=f-8\pi^2\sum_j\gamma_jG_0(q_j,\cdot)$, with the mean-zero normalization of \S\ref{sub:geometric}. Consider distributional conformal factors with $v-v_D^0\in H^2(M)$. Write $D_*=d_{\fp}(\kappa_*)$ as in \eqref{degree} and use the coefficients $a_q,\beta_q$ of \S\ref{sec:graded}.

\begin{theorem}[Existence and index for the geometric equations]
For either equation below, assume that the indicated mass is a positive wall and that \eqref{W} holds there.
\begin{align}
 \kappa_*=\kappa_D:\quad
 P_gv+2Q_g&=2he^{4v}-8\pi^2\sum_j\gamma_j\delta_{q_j},
 \label{geom-uncomp}\\
 \kappa_*=\kappa_P:\quad
 P_gv+2Q_g&=2he^{4v}-8\pi^2\sum_j\gamma_j(\delta_{q_j}-V^{-1}dV_g).
 \label{geom-comp}
\end{align}
In either case the degree in the mean-zero fixed-point formulation is $D_*$. There is a solution if $D_*\ne0$. At $D_*=0$, a solution exists if for some $j\ge0$
\begin{equation}\label{geom-zero-test}
 \sum_{q=0}^j(-1)^{j-q}(a_q-\beta_q)<0 .
\end{equation}
The regular remainders of the solutions form a compact set in $W^{4,s}(M)$ for every $1<s<s_*$. Each solution of \eqref{geom-uncomp} gives a metric $e^{2v}g$ with $Q$-curvature $h$ on $M\setminus\Sigma$ and order $\gamma_j$ at $q_j$.

If all genuine solutions are nondegenerate, the bounds of Corollaries~\ref{index-bounds78} and~\ref{count78} apply to the distinct conformal factors, using their reduced indices. A factor of reduced index $i$ has index $k+i+1$ for the unconstrained linearized form of \eqref{geom-uncomp} or \eqref{geom-comp} on $H^2(M)$.
\end{theorem}
\begin{proof}
Integrating either equation gives $\int_M he^{4v}\dV=\kappa_*$ in the indicated row. Since $K_D=he^{4v_D^0}$, the resulting bijection and its inverse are
\[
 v=u+v_D^0+\tfrac14\log\!\frac{\kappa_*}{\int_M K_De^{4u}\dV},
 \qquad u=v-v_D^0-\tfrac1V\int_M(v-v_D^0)\dV .
\]
Replacing $\kappa_*$ in the displayed lift by an arbitrary $\kappa>0$ defines $v_\kappa$. For $u$ solving \eqref{normalized} at that mass, the equations for $G_0$ and $f$ give
\begin{equation}
 P_gv_\kappa+2Q_g=2he^{4v_\kappa}
       -8\pi^2\sum_j\gamma_j\delta_{q_j}+2(\kappa_D-\kappa)V^{-1}dV_g .
\end{equation}
The constant background vanishes at $\kappa_D$ and is
$(8\pi^2/V)\sum_j\gamma_j$ at $\kappa_P$. The bijection transfers
compactness from Proposition~\ref{compact} because its normalization
constant is continuous on $W^{4,s}$. The Green expansion at $q_j$ gives $v=\gamma_j\log d_g(\cdot,q_j)+O(1)$, and conformal covariance gives the asserted curvature off the marks for \eqref{geom-uncomp}.

Theorem~\ref{deg} gives the degree and nonzero-degree existence. At degree zero the alternating sum in \eqref{geom-zero-test} is the coefficient of $t^j$ in $(A_--P_*)/(1+t)$. Apply Corollary~\ref{count78}. Neither argument requires invertibility at the genuine solutions.

Finally let $d\mu=he^{4v}dV_g/\kappa_*$. The derivative of the lift sends a volume-mean-zero variation $z$ to $z-\int z\,d\mu$. This is an isomorphism onto the $\mu$-mean-zero hyperplane, with inverse given by subtraction of the volume mean. On that hyperplane the form of $P_g-8he^{4v}$ pulls back to $D^2J_{\kappa_*}(u)/2$ by \eqref{hessian}. Its complementary constant line is orthogonal for this form and has value $-8\kappa_*c^2$ at the constant $c$. Thus the nullities agree, and \eqref{shift} gives the index assertion. The bijection preserves all counts.
\end{proof}

\subsection{A degree-zero example on the round four-sphere}
\begin{proposition}
On the unit round $S^4$, mark the two poles with orders $\gamma_1=\gamma_2=-\tfrac12$. There is a nonempty open set of smooth positive $h$ (in the $C^\infty$ topology) for which \eqref{geom-uncomp} has degree one, while \eqref{geom-comp} has degree zero and admits a solution. For coefficients in this set with nondegenerate genuine solutions, the compensated equation has a solution of reduced index one and a second solution of positive even reduced index.
\end{proposition}
\begin{proof}
Here $Q_g=3$, $V=8\pi^2/3$ and $P_g=\Delta_g^2-2\Delta_g$. Thus $k=0$, $\kappa_P=8\pi^2$ and $\kappa_D=4\pi^2$. The geometric equation lies on the first marked wall; the compensated equation lies on the shared wall of \S\ref{sec:graded}.

Let $F=\F_{1,\varnothing}$ and let $F^0$ be its value at $h=1$. On the round sphere, with $r=d_g(a,x)$, the normalized Green function is
\[
 G(a,x)=-2\log(2\sin(r/2))+\mathrm{const}.
\]
Indeed, in the variable $y=\cos r$ the radial Laplacian is $(1-y^2)\partial_y^2-4y\partial_y$. Applying $P_g$ to $-\log(2(1-y))$ gives $-6=-16\pi^2/V$ away from $a$, and the singular term $-2\log r$ gives $16\pi^2\delta_a$. The additive constant fixes the mean. The two pole contributions in \eqref{weight} now give
\[
 F^0(s)=-2\log(\cos s)+\mathrm{const},\qquad -\pi/2<s<\pi/2
\]
in latitude $s$, since $f=0$ and $H_g(a,a)$ is constant. Moreover $H_g(a,x)=\mathrm{const}+r^2/12+O(r^4)$, hence $4\Delta_xH_g(a,a)=8/3=2R_g/9$. By \eqref{allC}, at an ordinary critical point the sign of $C$ is exactly the sign of $\Delta_g F$.

Choose a smooth cutoff $\chi(s)$ supported away from the poles and equal to one near the equator. With $z\in S^3$ the angular variable on each latitude, prescribe
\[
 F(s,z)=F^0(s)+\epsilon\chi(s)(z_4-z_1^2).
\]
For sufficiently small $\epsilon>0$ its radial derivative has the sign of $s$ off the equator, since $(F^0)'=2\tan s$ and $\chi'$ is supported away from zero. On $S^3$, the function $z_4-z_1^2$ has two minima at $z_4=-1/2$, $z_1=\pm\sqrt3/2$, and critical points of index one and three at $-e_4$ and $e_4$. Indeed its constrained Hessian is the restriction of $\operatorname{diag}(-2,0,0,0)-2\lambda I$ where $\nabla(z_4-z_1^2)=2\lambda z$. The equatorial radial direction is positive. The ordinary interaction therefore has exactly four critical points, of indices $0,0,1,3$.

At the equator the radial eigenvalue is $2$ and the tangential eigenvalues are $\epsilon$ times the eigenvalues of the constrained Hessian on the unit $S^3$. The trace is positive at the minima and at the index-three point for sufficiently small $\epsilon$. At the index-one point the eigenvalues are $2,-\epsilon,\epsilon,\epsilon$. In normal coordinates with the first axis its negative direction, precompose $F$ with the time-one flow of a smooth cutoff of $(\log L)x_1\partial_{x_1}$ supported in a ball containing no other critical point. This diffeomorphism fixes the point and multiplies its negative eigenvalue by $L^2$ while preserving the others. Choosing $\epsilon L^2>2+2\epsilon$ makes its trace negative. Precomposition preserves all four critical points and their indices. The resulting $F-F^0$ vanishes in the caps and defines a smooth $\log h$ on the sphere.

Finally add $B r^2$ to $\log h$ at each pole, cut off in a ball of radius $\delta$ with a fixed rescaled cutoff equal to one near zero. This preserves the values and first jets there and adds $8B$ to the Laplacian at each pole. In \eqref{allC} the logarithmic coefficients of both single-mark configurations at the first wall and of the marked pair at the shared wall are positive linear combinations of these Laplacian increments plus fixed terms. First choose $B>0$ sufficiently large to make all three coefficients positive. Then choose $\delta$ small so that the derivative of the added term, bounded by $CB\delta$, is smaller than $|\partial_s F^0|$ in the punctured caps. No new ordinary critical point occurs and the four prescribed signs persist. Compactness of center configurations and continuation of Morse critical points in Proposition~\ref{compact} give openness of these data at both walls.

There are no retained ends at the first wall, so \eqref{degree} gives degree one. At the shared wall the only retained end is the ordinary critical point of index one. Consequently $A_-(t)=t^3$, while $P_*(t)=t$. Both Euler values are $-1$, and
\[
 \frac{A_-(t)-P_*(t)}{1+t}=t^2-t.
\]
The coefficient $-1$ of $t$ gives \eqref{geom-zero-test}. For Morse solutions, \eqref{strong78} at $j=1$ gives $m_1\ge m_0+1$. The Euler identity makes the total even and odd counts equal, hence a positive even index must also occur. Proposition~\ref{compact} gives a residual set of such Morse coefficients within the open set.
\end{proof}

\bigskip
\begingroup
\small\raggedright
\noindent\textbf{Mohameden Ahmedou}\\
Department of Mathematics, Justus Liebig University Giessen, Arndtstra\ss e 2, 35392 Gie\ss en, Germany.\\
\href{mailto:Mohameden.Ahmedou@math.uni-giessen.de}{Mohameden.Ahmedou@math.uni-giessen.de}\par\smallskip

\noindent\textbf{Lina Wu}\\
Department of Mathematics, Beijing Jiaotong University, Beijing 100044, China.\\
\href{mailto:lnwu@bjtu.edu.cn}{lnwu@bjtu.edu.cn}\par\smallskip

\noindent\textbf{Lei Zhang}\\
Department of Mathematics, University of Florida, 1400 Stadium Rd, Gainesville, FL 32611, USA.\\
\href{mailto:leizhang@ufl.edu}{leizhang@ufl.edu}\par
\endgroup

\end{document}